\documentclass[11pt]{article}
\usepackage[T1]{fontenc}
\usepackage{lmodern}
\usepackage{amsmath,amssymb,amsthm,mathtools}
\usepackage[margin=1in]{geometry}
\usepackage{microtype}
\usepackage{enumitem}
\usepackage{booktabs}
\usepackage{xcolor}
\usepackage[colorlinks=true,linkcolor=blue!50!black,citecolor=blue!50!black,urlcolor=blue!50!black]{hyperref}
\hypersetup{pdftitle={Elliptic Monotonicity Formulas for p-Green Functions}}
\numberwithin{equation}{section}
\newtheorem{theorem}{Theorem}[section]
\newtheorem{proposition}[theorem]{Proposition}
\newtheorem{lemma}[theorem]{Lemma}
\newtheorem{corollary}[theorem]{Corollary}
\theoremstyle{definition}

\newtheorem{assumption}[theorem]{Assumption}
\theoremstyle{remark}
\newtheorem{remark}[theorem]{Remark}

\newcommand{\Ric}{\operatorname{Ric}}
\newcommand{\Hess}{\nabla^2}
\newcommand{\Div}{\operatorname{div}}
\newcommand{\tr}{\operatorname{tr}}
\renewcommand{\Cap}{\operatorname{Cap}}
\newcommand{\Vol}{\operatorname{Vol}}

\newcommand{\dd}{\,d}
\newcommand{\A}{\mathsf A}
\newcommand{\D}{\mathcal D_p}
\newcommand{\E}{\mathcal E_{p,\beta}}
\newcommand{\LB}{\mathcal L_b}
\newcommand{\LG}{\mathcal L_G}
\newcommand{\LV}{\mathcal L_v}
\newcommand{\om}{\omega_{n-1}}
\newcommand{\bs}{\beta_*}
\newif\ifshowchanges
\showchangesfalse
\newcommand{\rev}[1]{\ifshowchanges
  {\color{red}\hypersetup{allcolors=red}#1}\else#1\fi}
\allowdisplaybreaks[2]
\setlist[enumerate]{itemsep=3pt}

\title{\bf Elliptic Monotonicity Formulas for $p$-Green Functions}
\author{Yu-Zhao Wang\footnote{Y.-Z. Wang was supported by the Fundamental Research Program of Shanxi Province (Grant No. 202303021211001).}}
\date{ }
\begin{document}
\maketitle
%\begin{center}
%\small\itshape
%Dedicated to the memory of my advisor, Professor Wenyi Chen
%\end{center}
\begin{abstract}
In this paper, we derive $p$-Laplace analogues of Colding’s three monotonicity formulas in the original metric on complete noncompact $n$-dimensional manifolds with nonnegative Ricci curvature, for $n\geq3$ and $1<p<n$. We assume a normalized positive minimal $p$-Green function vanishing at infinity with Euclidean pole asymptotics. Two complementary linearizations and the $p$-Bochner formula yield exact derivative identities with a common nonnegative defect involving an anisotropic trace-free Hessian and Ricci curvature. These identities give exterior cone rigidity and integral defect estimates, and recover the sharp gradient estimate and Euclidean rigidity of Mari–Rigoli–Setti. A one-parameter extension of Colding–Minicozzi type connects the formulas with known monotonicity and convexity identities for $p$-capacitary potentials.\vspace{3mm}


\vspace{3mm}


\end{abstract}

\section{Introduction}
Monotonicity formulas play an important role in geometric analysis. Colding's monotonicity formulas \cite{Colding} provide a connection
between harmonic Green functions and the geometry of manifolds
with nonnegative Ricci curvature. Their strength lies in the
derivative identities: the change of each monotone quantity is
expressed as an integral of nonnegative Hessian and curvature
terms. These identities characterize the equality case and
measure the deviation from cone geometry between different
scales. They therefore provide geometric information beyond
the sign of a derivative.

For a normalized positive Green function \(G\) in dimension
\(n\geq3\), set \(b=G^{1/(2-n)}\). This function agrees with distance
from the pole in Euclidean space. The identities
\begin{equation}\label{b2}
 \Delta b^2=2n|\nabla b|^2,\qquad
 |\nabla b^2|^2=4b^2|\nabla b|^2
\end{equation}
allow the Bochner formula to be organized around the defect
\begin{equation}\label{tracefree}
 \left|\nabla^2b^2-\frac{\Delta b^2}{n}g\right|^2
 +\Ric(\nabla b^2,\nabla b^2).
\end{equation}
The first term measures the failure of the cone Hessian identity,
while the second is nonnegative under the curvature assumption.
Colding--Minicozzi \cite{CM} extended this construction by varying
the gradient exponent in the level-set integrals. The role of
monotonicity and its derivative estimates in the study of
asymptotic geometry is further discussed in
\cite{CMTangent,CMMono}.

Parabolic entropy formulas provide a related motivation for
the present work. Perelman \cite{Perelman} introduced entropy
functionals for the Ricci flow and proved their monotonicity.
His \(\mathcal W\)-entropy formula expresses the derivative
as an integral of a geometric square, with shrinking gradient
Ricci solitons as the equality models. Inspired by this work,
Ni \cite{Ni} established an entropy formula for the linear
heat equation on a fixed Riemannian manifold. Its derivative
is the negative integral of a corrected Hessian square and
a Ricci curvature term, yielding monotonicity under
nonnegative Ricci curvature. Colding's elliptic formulas
were motivated in part by these parabolic identities and
their connections with sharp gradient estimates.

Kotschwar--Ni \cite{KN} extended the parabolic entropy theory
to the nonlinear equation
\[
 \frac{\partial}{\partial t}(u^{p-1})
 =(p-1)^{p-1}\Delta_pu,\qquad
 \Delta_pu:=\Div(|\nabla u|^{p-2}\nabla u).
\]
Using the linearized \(p\)-Laplacian and a nonlinear Bochner
identity, they obtained an entropy formula whose derivative
contains a Hessian square and a Ricci curvature term.
The Hessian norm is determined by the coefficient tensor
of the linearized operator. Their formula is monotone under
nonnegative Ricci curvature and recovers Ni's formula when
\(p=2\). Related gradient
estimates and entropy formulas for nonlinear diffusion
equations appear in \cite{LNVV,WC}; see also \cite{NiMo}
for a broader account of parabolic monotonicity.

Colding's elliptic formulas and Kotschwar--Ni's nonlinear
parabolic formula together motivate the present work.
Our aim is to extend Colding's three formulas to \(p\)-Green
functions while retaining their exact nonnegative geometric
defects. This question also arises naturally from nonlinear
potential theory. The \(p\)-harmonic equation is the
Euler--Lagrange equation of the \(p\)-Dirichlet energy, and
its conserved flux determines \(p\)-capacity. Monotonicity
formulas for its Green levels can therefore connect nonlinear
energy with area, mean curvature, and cone geometry. Keeping
the full derivative identities is essential for understanding
equality and estimating the geometric error between scales.

The main difficulty is to find quantities compatible with both
the nonlinear Bochner identity and the conserved \(p\)-harmonic
flux. For \(1<p<n\), the Euclidean \(p\)-Green function is proportional
to \(r^{-(n-p)/(p-1)}\), which suggests the radius
\[
 b=G^{-\frac{p-1}{n-p}}.
\]
The appropriate replacement for \(b^2\) is
\(v=b^{\frac{p}{p-1}}\). This power gives the nonlinear counterparts
of \eqref{b2} and allows the Bochner calculation to close. The linearization at \(v\) determines the Hessian norm.
Subtracting the trace with respect to its coefficient tensor
then produces a defect that vanishes on the Euclidean model.

A further issue appears when one tries to convert this
calculation into a product identity. Continuing with the
same operator on \(G|\nabla b|^p\) leaves additional drift
and zero-order terms. These terms disappear automatically
when \(p=2\), since the linearized operators all reduce to
the Laplacian. For general \(p\), matching the drift leads
instead to the linearization at \(G\). The identity
\(
 \mathcal L_GG=(p-1)\Delta_pG=0
\)
then removes the zero-order term. Thus the two linearizations
have complementary roles: the one at \(v\) produces the
Bochner defect, and the one at \(G\) gives a product identity
adapted to the Green function. The conserved flux supplies
the corresponding level-set measure and determines the
three monotone combinations.

This construction yields three exact formulas with a common
nonnegative defect and recovers Colding's formulas when
\(p=2\). 
Weighted Hessian estimates allow the identities to be
integrated across critical sets. Their consequences include
exterior cone rigidity, capacity and mean-curvature
comparisons, and annular integral estimates for the
trace-free Hessian and the other geometric defect terms.
After establishing this fixed-exponent theory, we vary the
gradient exponent in the spirit of Colding--Minicozzi.

In \cite{MRS}, Mari--Rigoli--Setti introduced a model-based fake distance
for the \(p\)-Green kernel and proved a sharp gradient bound
with model rigidity under Ricci lower bounds
\cite[Theorems~1.2 and~2.19]{MRS}. For the Euclidean model,
their fake distance is precisely our normalized Green radius
\(b\); their theorem therefore includes \(|\nabla b|\leq1\)
and its pointwise Euclidean equality case.
They also use a weighted Bochner identity with an
anisotropic trace-free decomposition
\cite[Propositions~2.14--2.15]{MRS}, and their Lemma~5.3
gives the corresponding flux averages and differentiation rule.

Related elliptic monotonicity results are already available
in nonlinear potential theory. Benatti--Fogagnolo--Mazzieri
\cite{BFM} proved monotonicity, convexity, and rigidity for
level-set integrals of \(p\)-capacitary potentials under
nonnegative Ricci curvature. Their approach uses a conformal
metric and a \(p\)-Bochner identity, and leads to sharp
geometric inequalities. Their later work \cite{BFMasym}
studies capacitary asymptotics on asymptotically conical
manifolds.
Benatti--Pluda--Pozzetta \cite{BPP} further developed the
regularity of the level sets and related these monotonicity
formulas to the weak inverse mean curvature flow.

In dimension three, Munteanu--Wang \cite{MW} established
Green-function comparisons under scalar-curvature assumptions.
Chan--Chu--Lee--Tsang \cite{CCLT} extended the associated
monotonicity to \(p\)-Green functions and proved rigidity
results. After a change of variable and normalization, their
functional agrees with the member \(\beta=3-p\) of our
parameter family. Under the stronger hypothesis
\(\Ric\geq0\), our identities express its first derivatives
as nonnegative volume-defect integrals.

For exterior potentials obtained from Green levels, the family
\(A_\beta\) agrees with the capacitary family in \cite{BFM}
after a change of variable and normalization. The third formula
corresponds to their convexity identity, and the other two follow
by weighted integration with the normalized pole limits.

The contribution of this paper is a direct, unified proof
in the original metric of Colding's three surface and
volume formulas. The calculation explains two choices:
the power \(v=b^{p/(p-1)}\), which makes the Bochner
identity close and selects \(B\), and the linearization
\(\LG\), which matches the drift and satisfies \(\LG G=0\).
The resulting product identity has no zero-order remainder
and is compatible with the conserved Green flux.
One local defect equation thus yields both the interior
formulas, whose constants are fixed at the pole, and the
exterior formula, whose boundary term is fixed by boundedness
at infinity. The proof uses no Hessian asymptotics at the
pole and no asymptotic expansion at infinity. Weighted
estimates justify the integrations across the critical set
and lead to rigidity and annular geometric-defect estimates.
The same product identity gives a direct proof of the
sharp theorem of Mari--Rigoli--Setti. Section~\ref{sec:comparison}
records the exact correspondences with \cite{MRS,BFM,CCLT};
the emphasis here is this common derivation and its
geometric consequences.

We now state the assumptions on the Green function and
the main results.

\begin{assumption}\label{ass:standing}
Let \((M^n,g)\) be smooth, connected, complete, and noncompact, with
\(n\geq3\) and \(\Ric\geq0\). Let \(1<p<n\), \(o\in M\), and
\(a=(n-p)/(p-1)\). Assume that \(G\) is a positive minimal
\(p\)-Green function normalized by
\[
 -\Delta_pG=a^{p-1}\om\delta_o,
 \qquad \om=|\mathbb S^{n-1}|.
\]
Assume that \(G\to0\) at infinity and that, uniformly near the pole,
\begin{equation}\label{eq:pole}
 b/d(o,\cdot)\longrightarrow1,\qquad
 |\nabla b|\longrightarrow1,\qquad b=G^{-1/a}.
\end{equation}
\end{assumption}
The assumed existence of a positive minimal \(p\)-Green
function already implies that \(M\) is \(p\)-nonparabolic
\cite{BFM}. The decay at infinity and the pole asymptotics
are explicit hypotheses. 
The normalized pole limits agree with
\cite[Proposition~2.13]{MRS}; no Hessian asymptotics at the
pole are assumed.
No separate nonparabolicity
assumption is needed. 
In particular, \(b\) extends
continuously by \(b(o)=0\) and is proper.
We do not assume Euclidean volume growth.
The Euclidean fundamental solution fixes the exponent \(a\) and
hence the radius \(b\). To state the Bochner identities, we use
\begin{equation}\label{eq:constants}
 q=\frac p{p-1},\qquad v=b^q,\qquad w=|\nabla b|,\qquad F=w^p.
\end{equation}
The power \(q\) is selected by the requirement that \(\Delta_pv\)
be a constant multiple of \(w^p\), with no remaining factor of
\(b\); this is verified in \eqref{eq:v-relations}.

On \(M\setminus(\{o\}\cup Z)\), where \(Z=\{w=0\}\), set
\[
 \nu=\frac{\nabla b}{w},\qquad
 \A^{ij}:=g^{ij}+(p-2)\nu^i\nu^j.
\]
We use the metric to identify \(\A\) with the corresponding
positive endomorphism when writing the operators below.
For \(f\in\{b,v,G\}\), define the linearization of
the \(p\)-Laplacian at \(f\) by
\begin{equation}\label{Linear}
 \mathcal L_f u
 :=\Div\bigl(|\nabla f|^{p-2}\A\nabla u\bigr).
\end{equation}
Thus
\[
 \left.\frac{d}{dt}\right|_{t=0}\Delta_p(f+tu)
 =\mathcal L_f u,
 \qquad
 \LB u=\Div(w^{p-2}\A\nabla u).
\]
The gradients of \(b,v,G\) are parallel, with that of \(G\)
pointing in the opposite direction. Since reversing the
gradient direction leaves \(\A\) unchanged, the same tensor
appears in all three linearizations.

The Hessian square in the linearized Bochner formula uses the following
trace and norm: for a symmetric tensor \(S\),
\[
 \tr_{\A}S:=\A^{ij}S_{ij},\qquad
 |S|_{\A}^2:=\A^{ik}\A^{j\ell}S_{ij}S_{k\ell},
 \qquad |\xi|_{\A}^2:=\A^{ij}\xi_i\xi_j
\]
for a covector \(\xi\). Accordingly, we subtract the
\(\A\)-trace, rather than the ordinary trace, and include the Ricci
term to define the nonnegative defect
\begin{equation}\label{Defect}
 \D:=|B|_{\A}^2+|\nabla v|^{2p-4}\Ric(\nabla v,\nabla v),
\end{equation}
where
\begin{equation}\label{eq:B-D}
 B:=|\nabla v|^{p-2}\Hess v-\frac{\Delta_pv}{n}\A^{-1}.
\end{equation}
This tensor satisfies \(\tr_{\A}B=0\). Here \(\A^{-1}\) is the covariant
inverse tensor, whose components are
\[
 (\A^{-1})_{ij}=g_{ij}+(q-2)\nu_i\nu_j,
 \qquad \A^{ik}(\A^{-1})_{kj}=\delta^i_j.
\]

These tensors are initially defined on \(w>0\). The scalar defect
is extended by zero on \(Z\) in integrals; this convention is
justified in Section~\ref{sec:critical}.

\begin{theorem}[Basic defect identities]
\label{thm:basic-main}
On \(w>0\), the following identities hold:
\begin{align}
 v\LV F+(p-1)(2-n)|\nabla v|^{p-2}
       \langle\nabla v,\nabla F\rangle
 &=\frac{p-1}{q^{p-1}}\D,
 \label{eq:first-basic}\\
 \LG(GF)
 &=(p-1)a^{p-2}q^{3-2p}b^{-n}\D,
 \label{eq:second-basic}\\
 \LB F-(p-1)(n+a-2)b^{-1}w^{p-2}
       \langle\nabla b,\nabla F\rangle
 &=(p-1)q^{3-2p}b^{-2}\D,
 \label{eq:fixed-master}\\
 \Div\bigl(b^{2-n-a}w^{p-2}\A\nabla F\bigr)
 &=(p-1)q^{3-2p}b^{-n-a}\D.
 \label{eq:fixed-div}
\end{align}
\end{theorem}
\begin{remark}The first two identities are the \(p\)-counterparts of
Colding's Lemma~2.7 in {\rm\cite{Colding}}.
The last two express the same defect equation in the
Green radius and in divergence form; they will be used
for the first and third monotonicity formulas.
A related weighted radius identity, with an anisotropic
trace-free decomposition, was obtained in
\cite[Proposition~2.15]{MRS}. Its precise correspondence
with our \(\beta=2\) identity is given in
Section~\ref{subsec:nonlinear-potential}.
The weak extensions needed for integration are established
in Section~\ref{sec:critical}.
\end{remark}

The normalized conserved flux on a level is
\(r^{1-n}w^{p-1}\,d\sigma\). Averaging the Bochner quantity
\(F=w^p\) against this flux determines the surface exponent
\(2p-1\); coarea then determines the matching volume exponent.
Thus the two quantities corresponding to Colding's \(A\) and \(V\)
are
\begin{equation}\label{eq:fixed-AV}
 A_p(r):=r^{1-n}\int_{b=r}w^{2p-1}\dd\sigma,
 \qquad
 V_p(r):=r^{-n}\int_{b<r}w^{2p}\dd\mu.
\end{equation}
The pole value of \(A_p\) is \(\om\). Subtracting this value and
taking a logarithmic primitive gives a useful alternative form of
the second formula:
\begin{equation}\label{eq:fixed-Vinf}
 V_{p,\infty}(r):=\int_1^r\frac{A_p(s)-\om}{s}\dd s.
\end{equation}
Changing the reference radius changes only an additive constant.
The surface integrals are initially defined for almost every
radius, over the noncritical part of each level. Statements at
every radius use the continuous representatives constructed in
Lemma~\ref{lem:fixed-averages}.

\begin{theorem}[Three Colding-type formulas]\label{thm:fixed}
Under Assumption~\ref{ass:standing}, the function
\(A_p\) has a
\(W^{2,1}_{\mathrm{loc}}(0,\infty)\), hence \(C^1\), representative.
For \(r>0\) and \(0<s<r\),
\begin{align}
 [A_p-q(n-1)V_p]'(r)
 &=q^{3-2p} r^{-n-1}\int_{b<r}\D\dd\mu,
 \label{eq:fixed-M1}\\
 [r^{-a}(A_p(r)-\om)]'
 &=q^{3-2p} r^{-a-1}\int_{b<r}b^{-n}\D\dd\mu,
 \label{eq:fixed-M2}\\
 [t^{1-a}A_p'(t)]_{t=s}^{t=r}
 &=q^{3-2p}\int_{s<b<r}b^{-n-a}\D\dd\mu.
 \label{eq:fixed-M3}
\end{align}
Thus the three quantities on the left are nondecreasing.
Equivalently, the second formula can be written as
\begin{equation}\label{eq:fixed-log}
 [A_p-aV_{p,\infty}]'(r)
 =\frac{q^{3-2p}}{r}\int_{b<r}b^{-n}\D\dd\mu.
\end{equation}
Moreover,
\begin{equation}\label{eq:fixed-exterior}
 A_p'(r)=-q^{3-2p}r^{a-1}\int_{b>r}b^{-n-a}\D\dd\mu,
\end{equation}
and \(A_p,V_p,V_{p,\infty}\) are nonincreasing. All displayed
integrals are finite. At the pole,
\(A_p\to\om\) and \(V_p\to\om/n\).
\end{theorem}

The following estimate and its pointwise equality case
are due to Mari--Rigoli--Setti
\cite[Theorems~1.2 and~2.19]{MRS}. We recover them from
\eqref{eq:second-basic}; the normalization is explained in
Section~\ref{sec:critical}.

\begin{theorem}[Sharp gradient bound (Mari--Rigoli--Setti \cite{MRS})]
\label{thm:sharp-main}
Under Assumption~\ref{ass:standing},
\[
 |\nabla b|\leq1\qquad\text{on }M\setminus\{o\}.
\]
Equality at one point away from \(o\) implies that \((M,g)\) is
Euclidean and \(G(x)=|x-o|^{-a}\).
\end{theorem}

\begin{remark}
When \(p=2\), \(A_p\) and \(V_p\) agree with
Colding's \(A\) and \(V\), and the three monotonicity
formulas recover those in \cite{Colding}.
The precise correspondence, including the logarithmic
and exterior forms, is recorded at the end of
Section~\ref{sec:fixed-formulas}.
\end{remark}

A concrete geometric consequence is averaged umbilicity
of the Green levels. Under Assumption~\ref{ass:standing},
if \(p<2+1/(n-2)\), then
\[
 R\int_R^{2R}\frac1{|\Sigma_t|}
       \int_{b=t}|\mathring h|^2\dd\sigma\dd t
 \longrightarrow0\qquad\text{as }R\to\infty.
\]
Here \(\Sigma_t=\{b=t\}\), and \(\mathring h\) is the
trace-free second fundamental form on the noncritical
part of the level; the integrals are understood for almost
every \(t\). In dimension three, this covers the full
range \(1<p<3\), without an assumption of Euclidean
volume growth. Theorem~\ref{thm:applications-main} gives
the annular estimate from which this conclusion follows.

Section~\ref{sec:beta} varies the gradient exponent in
the spirit of Colding--Minicozzi \cite{CM}. The exact
global defect identities and their rigidity statements
are proved for \(\beta>\beta_*=(n-p)/(n-1)\). At
\(\beta=\beta_*\), the limiting argument gives monotonicity
and convexity, without an exact absolutely continuous
defect representation. Section~\ref{sec:comparison}
discusses the geometric consequences and the relation
with capacitary, scalar curvature, and mass monotonicity.

Recently, Manea \cite{Manea} extended the sharp
gradient estimate and monotonicity formulas of
Colding and Colding--Minicozzi to closed manifolds
with a positive Ricci lower bound. His construction
uses the Green function of a curvature-dependent
Schr\"odinger operator and gives a sharp gradient
estimate with spherical rigidity. A natural question
is whether the \(p\)-Laplace framework developed here
admits analogous formulas under a positive Ricci
lower bound. This question is being investigated
in ongoing work by the author of this paper
\cite{WangPositiveRicci}.

Another natural direction is to extend the present results
to smooth metric measure spaces. In the linear case
\(p=2\), Song--Wei--Wu \cite{SWW} established
Colding-type monotonicity formulas under nonnegative
Bakry--\'Emery curvature assumptions. More recently,
the author and Z.-T.~Wu \cite{WW} extended Manea's
framework to closed weighted manifolds with a positive
finite-dimensional Bakry--\'Emery curvature lower bound.
These results, together with the nonlinear formulas
developed here, motivate the study of analogous
monotonicity formulas and rigidity for the weighted
\(p\)-Laplacian.

\section{Linearized Bochner identities and the geometry of the level sets}
\label{sec:operators}

\subsection{Linearized Bochner identities }
All calculations in this section take place on \(w>0\).
We first explain the choice of \(v=b^q\) and the
anisotropic trace-free tensor. We then prove
Theorem~\ref{thm:basic-main} and identify the geometry
measured by its defect.
Linearized Bochner formulas and anisotropic trace
decompositions are already used in \cite{KN,MRS,FMP,BFM}.
Here the choice of \(v\), the norm \(|\cdot|_{\A}\), and the
change to \(\LG\) are organized to produce the three
Colding-type formulas from one tensor.

Since \(G=b^{-a}\) and \((a+1)(p-1)=n-1\),
\[
 |\nabla G|^{p-2}\nabla G
 =-a^{p-1}b^{1-n}w^{p-2}\nabla b.
\]
The equation \(\Delta_pG=0\) therefore gives
\begin{equation}\label{eq:div-X}
 \Div(b^{1-n}w^{p-2}\nabla b)=0.
\end{equation}
Expanding the derivative of \(b^{1-n}\),
\[
 0=b^{1-n}\Delta_pb+(1-n)b^{-n}w^p,
\]
so, on the regular set,
\begin{equation}\label{eq:b-equation}
 \Delta_pb=(n-1)b^{-1}w^p.
\end{equation}
For \(\alpha>0\), the chain rule now gives
\begin{align*}
 \Delta_p(b^\alpha)
 &=\alpha^{p-1}b^{(\alpha-1)(p-1)}
       \left[\Delta_pb+(\alpha-1)(p-1)b^{-1}F\right]\\
 &=\alpha^{p-1}\bigl[n-p+(p-1)\alpha\bigr]
       b^{(p-1)\alpha-p}F.
\end{align*}
The exponent of \(b\) is zero precisely when \(\alpha=q\). Hence, for $v=b^q$, $F=w^p=|\nabla b|^p$,
\begin{equation}\label{eq:v-relations}
 \Delta_pv=nq^{p-1}F,\qquad
 |\nabla v|^p=q^pvF.
\end{equation}
These two relations are the nonlinear substitutes for the two
identities for \(b^2\) in \eqref{b2}. They identify both the
power \(v\) and the gradient quantity \(F\) before the Bochner
formula is applied.

To determine the Hessian norm in that formula, we first compute
the linearization. For \(f\in\{b,v,G\}\), differentiation of the flux
with respect to a variation \(f+t\phi\) gives
\begin{align*}
 \left.\frac d{dt}\right|_{t=0}
       |\nabla(f+t\phi)|^{p-2}\nabla(f+t\phi)
 =|\nabla f|^{p-2}\nabla\phi
 +(p-2)|\nabla f|^{p-4}
       \langle\nabla f,\nabla\phi\rangle\nabla f.
\end{align*}
Taking the divergence gives \(\mathcal L_f\phi\) defined in \eqref{Linear}. The operator is
linear in \(\phi\), although its coefficients depend nonlinearly on
\(f\). In a frame with first vector parallel to \(\nabla f\),
\(\A=\operatorname{diag}(p-1,1,\ldots,1)\). Therefore
\[
 \A\nabla f=(p-1)\nabla f,\quad
 \mathcal L_ff=(p-1)\Delta_pf,\quad
 \Delta_pf=|\nabla f|^{p-2}\tr_{\A}\Hess f.
\]
Thus the natural trace of
\(|\nabla v|^{p-2}\Hess v\) is \(\Delta_pv\).
Since
\[
 |\A^{-1}|_{\A}^2=n,\qquad
 \left\langle|\nabla v|^{p-2}\Hess v,
                  \A^{-1}\right\rangle_{\A}=\Delta_pv,
\]
the tensor \(B\) defined in \eqref{eq:B-D} satisfies
\begin{equation}\label{eq:B-norm}
 \tr_{\A}B=0,\qquad
 |B|_{\A}^2
 =|\nabla v|^{2p-4}|\Hess v|_{\A}^2
       -\frac{(\Delta_pv)^2}{n}.
\end{equation}
The Euclidean model explains this subtraction. When \(b=r\),
\[
 \Hess r^q=qr^{q-2}\bigl(g+(q-2)\,dr\otimes dr\bigr)
          =qr^{q-2}\A^{-1},
\]
\rev{It follows that}
\[
 |\nabla r^q|^{p-2}\Hess(r^q)=q^{p-1}\A^{-1},
 \qquad
 \Delta_p(r^q)=nq^{p-1},
\]
because \(q-1=1/(p-1)\). Thus
 \[B=|\nabla r^q|^{p-2}\Hess r^q-\frac{\Delta_pr^q}n\A^{-1}=0.\]
 The ordinary trace-free
Hessian would not vanish on this model when \(p\ne2\). This is the reason for using
the \(\A\)-trace in the defect.

\begin{proof}[\bf Proof of Theorem \ref{thm:basic-main}]
Write \(s=|\nabla v|\) during the Bochner calculation.
Covariant differentiation and commutation give
\[
 \nabla_k(\Delta_pv)
 =\nabla_i\bigl(s^{p-2}\A^{ij}v_{jk}\bigr)
       -s^{p-2}\Ric_{k\ell}v_\ell.
\]
Multiplication by \(v_k\) and the product rule yield
\[
 \frac12\LV(s^2)
 =\langle\nabla v,\nabla\Delta_pv\rangle
  +s^{p-2}\left[
       |\Hess v|^2+(p-2)|\nabla s|^2
       +\Ric(\nabla v,\nabla v)\right].
\]
The scalar chain rule gives
\[
 \frac1p\LV(s^p)
 =\frac12s^{p-2}\LV(s^2)
       +(p-2)s^{2p-4}|\nabla s|_{\A}^2.
\]
In a frame with first vector \(\nu\), we have
\(\A=\operatorname{diag}(p-1,1,\ldots,1)\) and
\(\nabla_i s=v_{1i}\). Therefore
\[
 |\Hess v|_{\A}^2
 =|\Hess v|^2
       +(p-2)\bigl(|\nabla s|^2+|\nabla s|_{\A}^2\bigr).
\]
Combining these identities gives the $p$-Bochner formula
(cf.~\cite{KN})
\begin{equation}\label{eq:p-Bochner}
 \frac1p\LV(|\nabla v|^p)
 =|\nabla v|^{2p-4}
       \left[|\Hess v|_{\A}^2
             +\Ric(\nabla v,\nabla v)\right]
  +|\nabla v|^{p-2}
       \langle\nabla v,\nabla\Delta_pv\rangle.
\end{equation}
 Substitute
\eqref{eq:B-norm} and \eqref{eq:v-relations} into the right-hand
side of \eqref{eq:p-Bochner}:
\begin{equation}\label{eq:bochner-expanded}
 \frac1p\LV(|\nabla v|^p)
 =\D+nq^{2p-2}F^2+nq^{p-1}|\nabla v|^{p-2}\langle\nabla v,\nabla F\rangle.
\end{equation}
For any two functions \(f,h\), the linear divergence operator has
the product rule
\[
 \LV(fh)=f\LV h+h\LV f+2|\nabla v|^{p-2}\langle\A\nabla f,\nabla h\rangle.
\]
Apply it to \(|\nabla v|^p=q^pvF\). Since
\(\LV v=(p-1)nq^{p-1}F\), we obtain
\begin{align}\label{LVF}
 \frac1p\LV(|\nabla v|^p)= \frac{q^p}p\LV(vF)
 &=\frac{q^p}{p}\bigl[v\LV F+F\LV v
                 +2(p-1)|\nabla v|^{p-2}\langle\nabla v,\nabla F\rangle\bigr]\notag\\
 &=\frac{q^{p-1}}{p-1}v\LV F+nq^{2p-2}F^2
       +2q^{p-1}|\nabla v|^{p-2}\langle\nabla v,\nabla F\rangle.
\end{align}
Combining \eqref{eq:bochner-expanded} with \eqref{LVF}, we obtain
\eqref{eq:first-basic}.

We next explain the change of operator. Since \(G=v^{-a/q}\),
the composition rule gives
\begin{align*}
 \LV G
 &=-\frac aq\frac Gv\LV v
   +\frac aq\left(\frac aq+1\right)\frac G{v^2}
          |\nabla v|^{p-2}|\nabla v|_{\A}^2\\
 &=(p-1)a q^{p-2}\frac{GF}{v}\,[a+q-n]\\
 &=na q^{p-2}(2-p)\frac{GF}{v},
\end{align*}
where \(a+q=n/(p-1)\). Thus \(G\) is generally not
\(\LV\)-harmonic. Expanding the product and inserting
\eqref{eq:first-basic} gives
\begin{equation}\label{eq:wrong-operator}
\begin{split}
 \LV(GF)={}&\frac{p-1}{q^{p-1}}\frac Gv\D
 +\frac{n(p-1)(p-2)}p\frac{G}{v}|\nabla v|^{p-2}
       \langle\nabla v,\nabla F\rangle
   +na q^{p-2}(2-p)\frac Gv F^2.
\end{split}
\end{equation}
Both extra terms vanish when \(p=2\). For general \(p\),
the drift has no fixed sign, so this product does not
give the required nonnegative defect identity.

To absorb the drift into the coefficient, multiply
the coefficient of \(\LV\) by a positive function
\(\theta(v)\), and define
\(L^{\theta}=\Div(\theta |\nabla v|^{p-2}\A\nabla\cdot)\).
Its product rule reads
\begin{equation}\label{eq:theta-product}
 L^{\theta}(GF)=\theta G\bigl[\LV F+(p-1)|\nabla v|^{p-2}
 \langle\nabla\log\theta+2\nabla\log G,\nabla F\rangle\bigr]
       +F L^{\theta}G.
\end{equation}
To match the drift in \eqref{eq:first-basic}, require
\[
 \nabla\log\theta+2\nabla\log G=(2-n)\nabla\log v.
\]
Since \(\nabla\log G=-(a/q)\nabla\log v\), this becomes
\[
 \nabla\log\theta=-\frac{n(p-2)}p\nabla\log v.
\]
Thus \(\theta=Cv^{-n(p-2)/p}\). The actual gradient ratio is
\[
 \frac{|\nabla G|}{|\nabla v|}
 =\frac aq b^{-(a+q)}
 =\frac aq b^{-n/(p-1)},
\]
so choosing \(C=(a/q)^{p-2}\) gives
\begin{equation}\label{eq:theta}
 \theta=\left(\frac aq\right)^{p-2}b^{-n(p-2)/(p-1)},
 \qquad \theta |\nabla v|^{p-2}=|\nabla G|^{p-2}.
\end{equation}
The required operator is therefore \(\LG\).
It also satisfies \(\LG G=(p-1)\Delta_pG=0\), eliminating the last term in
\eqref{eq:theta-product}.
For clarity, the two exponent cancellations are
\begin{align*}
 -\frac{n(p-2)}{p-1}-2a
 &=\frac{p(2-n)}{p-1}=q(2-n),\\
 -\frac{n(p-2)}{p-1}-a-q&=-n.
\end{align*}
Consequently \(\theta G/v=(a/q)^{p-2}b^{-n}\).
Substitution of \eqref{eq:first-basic} into
\eqref{eq:theta-product} yields
\[
 \LG(GF)
 =\frac{p-1}{q^{p-1}}\frac{\theta G}{v}\D
 =(p-1)a^{p-2}q^{3-2p}b^{-n}\D,
\]
which proves \eqref{eq:second-basic}. Thus the same choice
matches the drift and removes the zero-order term.

For the radius form, use
\[
 |\nabla v|^{p-2}=q^{p-2}b^{2-q}w^{p-2}.
\]
The product rule for the radial coefficient gives
\[
 \LV F=q^{p-2}b^{2-q}
 \left[\LB F+(p-1)(2-q)b^{-1}w^{p-2}
                      \langle\nabla b,\nabla F\rangle\right].
\]
Substitute this into \eqref{eq:first-basic}, use
\(\nabla v=qb^{q-1}\nabla b\), and divide by
\(q^{p-2}b^2\). Since
\[
 (2-q)+q(2-n)=2-q(n-1)=-(n+a-2),
\]
we obtain \eqref{eq:fixed-master}.

Finally,
\begin{align*}
 \Div\bigl(b^{2-n-a}w^{p-2}\A\nabla F\bigr)
 =b^{2-n-a}\left[
       \LB F-(p-1)(n+a-2)b^{-1}w^{p-2}
             \langle\nabla b,\nabla F\rangle\right].
\end{align*}
Equation~\eqref{eq:fixed-master} now gives
\eqref{eq:fixed-div}.
\end{proof}
\subsection{The geometry of the level sets}
We finish by expressing the defect in terms of the
geometry of the level sets. This identifies the quantities
controlled by the monotonicity formulas and used in
the rigidity arguments.

\begin{proposition}[Geometric form of the defect]
\label{prop:defect-geometry}
On a regular level, let \(h\) be the second fundamental
form with normal \(\nu\), let \(H=\tr h\), and set
\[
 \mathring h=h-\frac{H}{n-1}g_\Sigma.
\]
Then
\begin{equation}\label{eq:D-geometric}
 \D=q^{2p-2}b^2w^{2p-2}\left[
       |\mathring h|^2+\Ric(\nu,\nu)
       +2(p-1)|\nabla^\Sigma\log w|^2
       +\frac{n}{n-1}
          \left(H-\frac{n-1}{b}w\right)^2
       \right].
\end{equation}
\end{proposition}

\proof

We now compute the components and norm of \(B\).
At a point where \(w>0\), choose an orthonormal frame
\(e_1=\nu\), with \(e_i\), \(2\leq i\leq n\), tangent
to the level set of \(b\). We use the convention
\[
 h_{ij}=\langle\nabla_{e_i}\nu,e_j\rangle,
 \qquad H=\sum_{i=2}^n h_{ii}.
\]
All norms of tangential tensors below are taken with respect
to the induced metric \(g_\Sigma\).

Since \(w^2=|\nabla b|^2\), differentiation in any direction
\(X\) gives
\[
 2wX(w)
 =X(|\nabla b|^2)
 =2\Hess b(X,\nabla b)
 =2w\Hess b(X,\nu).
\]
Thus
\[
 \Hess b(\nu,\nu)=\partial_\nu w,
 \qquad
 \Hess b(\nu,e_i)=\Hess b(e_i,\nu)
 =\nabla_i^\Sigma w.
\]
For the tangential components, differentiating
\(\nabla b=w\nu\) yields
\begin{align*}
 \Hess b(e_i,e_j)
 =\langle\nabla_{e_i}(w\nu),e_j\rangle
 =e_i(w)\langle\nu,e_j\rangle
   +w\langle\nabla_{e_i}\nu,e_j\rangle
 =wh_{ij}.
\end{align*}
Consequently,
\begin{equation}\label{Db}
 \Delta b=\partial_\nu w+wH.
\end{equation}

We also have
\[
 db(\nu)=w,\qquad db(e_i)=0,
 \qquad
 g(\nu,\nu)=1,\qquad g(\nu,e_i)=0.
\]
Differentiating \(v=b^q\) gives
\[
 |\nabla v|^{p-2}\Hess v
 =q^{p-1}w^{p-2}
       \left[b\Hess b+(q-1)db\otimes db\right].
\]
Subtracting
\[
 \frac{\Delta_pv}{n}\A^{-1}
 =q^{p-1}w^p\bigl[g+(q-2)\nu\otimes\nu\bigr]
\]
and using \(db\otimes db=w^2\nu\otimes\nu\), we obtain
\begin{equation}\label{eq:B-T}
 B=q^{p-1}bw^{p-2}
       \left(\Hess b-\frac{w^2}{b}g
                    +\frac1b db\otimes db\right).
\end{equation}
Substituting these relations into \eqref{eq:B-T}, we obtain
\begin{align*}
 B_{\nu\nu}
 &=q^{p-1}bw^{p-2}
   \left(\partial_\nu w-\frac{w^2}{b}
                         +\frac{w^2}{b}\right)=q^{p-1}bw^{p-2}\partial_\nu w,\\
 B_{\nu i}
 &=q^{p-1}bw^{p-2}\nabla_i^\Sigma w,\\
 B_{ij}
 &=q^{p-1}bw^{p-2}
   \left(wh_{ij}-\frac{w^2}{b}(g_\Sigma)_{ij}\right).
\end{align*}
In particular,
\begin{equation}\label{Bnu}
 B(\nu,\cdot)=q^{p-1}bw^{p-2}\,dw.
\end{equation}
The equation for \(b\) relates the normal derivative of \(w\)
to the mean curvature. Indeed, by \eqref{Db},
\begin{align*}
 \Delta_pb
 &=\Div(w^{p-2}\nabla b)\\
 &=w^{p-2}\Delta b
   +(p-2)w^{p-3}\langle\nabla w,\nabla b\rangle\\
 &=w^{p-2}
   \left[\partial_\nu w+wH+(p-2)\partial_\nu w\right]\\
 &=w^{p-2}\left[(p-1)\partial_\nu w+wH\right].
\end{align*}
Together with \eqref{eq:b-equation}, this gives
\begin{equation}\label{eq:normal-relation}
 (p-1)\partial_\nu w+wH=(n-1)\frac{w^2}{b}.
\end{equation}
Equivalently,
\[
 (p-1)\partial_\nu w
 =-w\left(H-\frac{n-1}{b}w\right).
\]
This relation also verifies the anisotropic trace-free
property directly:
\[
 \tr_{\A}B
 =q^{p-1}bw^{p-2}
   \left[(p-1)\partial_\nu w+wH
                        -(n-1)\frac{w^2}{b}\right]
 =0.
\]

To compute the norm, observe that in the chosen frame
\[
 \A=\operatorname{diag}(p-1,1,\ldots,1).
\]
The definition
\[
 |B|_{\A}^2=\A^{\alpha\gamma}\A^{\beta\delta}
                         B_{\alpha\beta}B_{\gamma\delta}
\]
therefore gives
\[
 |B|_{\A}^2
 =(p-1)^2B_{\nu\nu}^2
 +2(p-1)\sum_{i=2}^n B_{\nu i}^2
 +\sum_{i,j=2}^n B_{ij}^2.
\]
The factor \(2\) in the mixed term accounts for both
\(B_{\nu i}\) and \(B_{i\nu}\). Substitution of the component
formulas yields
\begin{align*}
 \frac{|B|_{\A}^2}{q^{2p-2}b^2w^{2p-4}}
 ={}&(p-1)^2(\partial_\nu w)^2
     +2(p-1)|\nabla^\Sigma w|^2+\sum_{i,j=2}^n
       \left(wh_{ij}-\frac{w^2}{b}(g_\Sigma)_{ij}\right)^2.
\end{align*}
We next separate the tangential trace. Since
\[
 h_{ij}=\mathring{h}_{ij}
                +\frac{H}{n-1}(g_\Sigma)_{ij},
\]
we have
\[
 wh_{ij}-\frac{w^2}{b}(g_\Sigma)_{ij}
 =w\mathring{h}_{ij}
   +\left(\frac{wH}{n-1}-\frac{w^2}{b}\right)(g_\Sigma)_{ij}.
\]
The two terms are orthogonal because
\[
 \sum_{i,j=2}^n\mathring{h}_{ij}(g_\Sigma)_{ij}=0,
 \qquad
 \sum_{i,j=2}^n(g_\Sigma)_{ij}^2=n-1.
\]
Hence
\begin{align*}
 \sum_{i,j=2}^n
       \left(wh_{ij}-\frac{w^2}{b}(g_\Sigma)_{ij}\right)^2
 =&w^2|\mathring{h}|^2
   +\frac{w^2}{n-1}
          \left(H-\frac{n-1}{b}w\right)^2.
\end{align*}
Combining this with the norm calculation gives
\begin{align*}
 \frac{|B|_{\A}^2}{q^{2p-2}b^2w^{2p-4}}
 =(p-1)^2(\partial_\nu w)^2
     +2(p-1)|\nabla^\Sigma w|^2
     +w^2|\mathring{h}|^2+\frac{w^2}{n-1}
          \left(H-\frac{n-1}{b}w\right)^2.
\end{align*}
By \eqref{eq:normal-relation},
\[
 (p-1)^2(\partial_\nu w)^2
 =w^2\left(H-\frac{n-1}{b}w\right)^2.
\]
Moreover,
\[
 |\nabla^\Sigma w|^2
 =w^2|\nabla^\Sigma\log w|^2.
\]
Thus
\[
 |B|_{\A}^2
 =q^{2p-2}b^2w^{2p-2}
 \left[
 |\mathring{h}|^2
 +2(p-1)|\nabla^\Sigma\log w|^2
 +\frac{n}{n-1}
       \left(H-\frac{n-1}{b}w\right)^2
 \right].
\]
The coefficient \(n/(n-1)\) is the sum of the normal
contribution \(1\) and the tangential trace contribution
\(1/(n-1)\).

It remains to compute the curvature term. Since
\(
 \nabla v=qb^{q-1}w\nu,
\)
we have
\begin{align*}
 |\nabla v|^{2p-4}\Ric(\nabla v,\nabla v)
 =|\nabla v|^{2p-2}\Ric(\nu,\nu)
 =q^{2p-2}b^2w^{2p-2}\Ric(\nu,\nu).
\end{align*}
Adding this term to the expression for \(|B|_{\A}^2\)
proves \eqref{eq:D-geometric}.

The four terms in \eqref{eq:D-geometric} measure
umbilicity, radial Ricci curvature, tangential variation
of the gradient, and the cone mean-curvature relation.
Under \(\Ric\geq0\), each is nonnegative.
\endproof

\section{Sharp gradient control and weak identities}
\label{sec:critical}

We first give another proof of the sharp Green-radius
estimate and pointwise rigidity of Mari--Rigoli--Setti
\cite[Theorems~1.2 and~2.19]{MRS}. Our argument uses the
two facts supplied by the product identity
\eqref{eq:second-basic}:
\[
 \LG(GF)\geq0,\qquad \LG G=(p-1)\Delta_pG=0
 \qquad\text{on }w>0.
\]
Comparing \(GF\) with a constant multiple of \(G\)
recovers the sharp constant. The maximum principle is used
only on compact regions where \(w>0\), so no extension of
the singular linearized operator across the critical set
is needed for this proof.

For the monotonicity formulas, a further argument is needed.
Their divergence identities must hold across the critical
set \(Z=\{w=0\}\). We obtain a weighted Hessian estimate
by regularization and then use it to remove a cutoff
around \(Z\).

\subsection{The sharp gradient estimate}

To compare normalizations, write \(\mathcal G\) for the
minimal kernel in \cite{MRS}, normalized by
\(-\Delta_p\mathcal G=\delta_o\). Then
\[
 G=a\om^{1/(p-1)}\mathcal G,
 \qquad
 \mathcal G^{\mathbb R^n}(t)
 =\frac{t^{-a}}{a\om^{1/(p-1)}}.
\]
Their Euclidean-model fake distance \(\varrho\) is defined by
\(\mathcal G=\mathcal G^{\mathbb R^n}(\varrho)\), hence
\(\varrho=G^{-1/a}=b\). Their theorem includes both the
estimate below and its pointwise equality case, and also
allows more general Ricci lower bounds without assuming
\(G\to0\) at infinity. The argument here uses that decay
assumption to give a direct comparison proof adapted to
the same \(\LG\) product identity as our monotonicity
formulas. We use \cite{MRS} with its 2023 corrigendum,
which leaves this sharp Green-radius theorem unchanged.

We first establish a global bound for \(F\). Its purpose
is to control \(GF\) at infinity; the constant in this
preliminary estimate will not enter the sharp bound.

\begin{lemma}[A preliminary global bound]
\label{lem:rough}
Under Assumption~\ref{ass:standing},
\[
 b(x)\leq d(o,x),\qquad
 |\nabla b|(x)\leq C(n,p)
 \qquad(x\ne o).
\]
\end{lemma}

\begin{proof}
The first comparison is also \cite[Proposition~2.10]{BFM};
we include its proof to fix the normalization.
Write \(\rho=d(o,\cdot)\). Away from the pole and the
cut locus,
\[
 |\nabla(\rho^{-a})|^{p-2}\nabla(\rho^{-a})
 =-a^{p-1}\rho^{1-n}\nabla\rho.
\]
Consequently,
\[
 \Delta_p(\rho^{-a})
 =-a^{p-1}\rho^{1-n}
       \left(\Delta\rho-\frac{n-1}{\rho}\right).
\]
Since \(\Ric\geq0\), Laplace comparison gives
\[
 \Delta\rho\leq\frac{n-1}{\rho}.
\]
The comparison holds weakly across the cut locus, so
\(\rho^{-a}\) is \(p\)-subharmonic away from \(o\).

Fix \(0<\varepsilon<1\) and \(R>0\). On a punctured
geodesic ball, compare \(G\) with
\[
 (1-\varepsilon)(\rho^{-a}-R^{-a}).
\]
This function vanishes on \(\rho=R\), where \(G>0\).
The pole normalization gives
\[
 \frac{G}{\rho^{-a}}
 =\left(\frac{\rho}{b}\right)^a\longrightarrow1,
\]
so the comparison function also lies below \(G\) on
every sufficiently small inner sphere. The weak comparison
principle therefore gives
\[
 G(x)\geq
 (1-\varepsilon)(\rho(x)^{-a}-R^{-a})
 \qquad(0<\rho(x)<R).
\]
Letting \(R\to\infty\) and then \(\varepsilon\downarrow0\),
we obtain \(G\geq\rho^{-a}\), or equivalently \(b\leq\rho\).

For \(x\ne o\), the ball
\(B(x,\rho(x)/2)\) does not contain the pole.
The local gradient estimate of Wang--Zhang \cite{WZ}
applied to \(G\) on this ball gives
\[
 \frac{|\nabla G|(x)}{G(x)}
 \leq\frac{C(n,p)}{\rho(x)}.
\]
Differentiating \(b=G^{-1/a}\), we have
\[
 w=\frac{b}{a}\frac{|\nabla G|}{G}.
\]
Using \(b\leq\rho\) proves the claimed global bound.
\end{proof}

The choice of \(\LG\) is useful here for the same reason
as in the product identity: it annihilates \(G\).
Thus subtracting a constant multiple of \(G\) from \(GF\)
preserves the differential inequality. This gives a direct
comparison proof of Theorem~\ref{thm:sharp-main}.

\begin{proof}[\bf A direct proof of Theorem~\ref{thm:sharp-main}]
Fix \(\delta>0\), and define
\[
 H_\delta=GF-(1+\delta)G
 =G(F-1-\delta).
\]
On \(w>0\), linearity of \(\LG\) and
\eqref{eq:second-basic} give
\[
 \LG H_\delta
 =(p-1)a^{p-2}q^{3-2p}b^{-n}\D
 \geq0.
\]

We first check the behavior at the two ends of the manifold.
Near the pole, \(F\to1\). Hence, sufficiently close to \(o\),
\[
 F-1-\delta\leq-\frac{\delta}{2}.
\]
Since \(G\to\infty\) there,
\[
 H_\delta\longrightarrow-\infty
 \qquad\text{as }x\to o.
\]
At infinity, Lemma~\ref{lem:rough} bounds \(F\), while
\(G\to0\). Therefore
\[
 H_\delta\longrightarrow0
 \qquad\text{as }x\to\infty.
\]

Suppose that \(H_\delta(x)>0\) at some point \(x\).
Choose \(c\) with
\[
 0<c<H_\delta(x),
\]
and let \(U\) be the connected component of
\(\{H_\delta>c\}\) containing \(x\).
The preceding limits show that \(\overline U\) is compact
and does not contain the pole.

On \(\overline U\), we have
\[
 F=1+\delta+\frac{H_\delta}{G}\geq1+\delta.
\]
In particular \(w>0\) throughout a neighborhood of
\(\overline U\). Since
\[
 |\nabla G|=ab^{-a-1}w,
\]
the gradient of \(G\) is bounded away from zero on this
compact set. The eigenvalues of \(\A\) are \(1\) and
\(p-1\), so the coefficients of \(\LG\) are smooth and
uniformly elliptic there.

We may therefore apply the maximum principle to
\(\LG H_\delta\geq0\) on \(U\). Since
\(H_\delta=c\) on \(\partial U\), it gives
\(H_\delta\leq c\) in \(U\), contradicting the choice of \(x\).
Thus \(H_\delta\leq0\) everywhere, and positivity of \(G\)
gives
\[
 F\leq1+\delta.
\]
Letting \(\delta\downarrow0\), we obtain
\[
 F=w^p\leq1.
\]
Hence \(|\nabla b|\leq1\).

We now prove the equality case. Suppose that \(F(x_0)=1\)
at some \(x_0\ne o\). The function
\[
 H_0=G(F-1)
\]
is nonpositive and vanishes at \(x_0\). On the connected
component \(U\) of \(\{w>0\}\) containing \(x_0\), it satisfies
\[
 \LG H_0\geq0.
\]
The strong maximum principle gives \(H_0=0\) on \(U\),
and therefore \(w=1\) there.

Continuity makes \(U\) closed in \(M\setminus\{o\}\).
Indeed, a limit point of \(U\) still has \(w=1\), and a
small connected neighborhood of that point lies in
\(\{w>0\}\) and belongs to the same component.
Since \(M\setminus\{o\}\) is connected, we conclude that
\[
 w=1\qquad\text{on }M\setminus\{o\}.
\]
Thus \(H_0=0\) everywhere away from the pole.
The product identity then gives \(\D=0\), and hence \(B=0\).
Equation~\eqref{eq:B-T} becomes
\[
 b\Hess b=g-db\otimes db.
\]

We next identify \(b\) with the distance from \(o\).
Starting at \(x\ne o\), follow the flow of \(-\nabla b\).
Along this curve,
\[
 \frac{d}{dt}b=-|\nabla b|^2=-1,
 \qquad |\dot\gamma|=1.
\]
Properness of \(b\) permits the flow to continue through
every band with positive lower endpoint. Thus it exists
for \(0\leq t<b(x)\), with
\[
 b(\gamma(t))=b(x)-t.
\]
The curve has finite length \(b(x)\). Completeness gives
a limit as \(t\uparrow b(x)\), and continuity of \(b\)
makes its value at that limit zero. Since \(b=0\) only
at \(o\), the limit is \(o\). Therefore
\[
 d(o,x)\leq b(x).
\]
Together with Lemma~\ref{lem:rough}, this proves
\(b=d(o,\cdot)\).

Set \(r=b\). The gradient flow identifies the positive
levels, and the metric takes the form
\[
 g=dr^2+g_r.
\]
For vectors tangent to a level and transported by the flow,
the Hessian identity gives
\[
 \frac{d}{dr}g_r
 =2\Hess b|_{T\Sigma_r}
 =\frac2r g_r.
\]
Hence \(r^{-2}g_r\) is independent of \(r\).
Near the smooth pole, geodesic polar coordinates give
\[
 r^{-2}g_r\longrightarrow g_{\mathbb S^{n-1}},
\]
where \(g_{\mathbb S^{n-1}}\) is the unit round metric.
It follows that
\[
 g=dr^2+r^2g_{\mathbb S^{n-1}}
\]
on all of \(M\setminus\{o\}\).
Thus \(M\) is Euclidean, and
\[
 G=b^{-a}=|x-o|^{-a}.
\]
\end{proof}

\begin{remark}[Sharpness and the comparison function]
\label{rem:sharp}
The estimate is equivalently
\[
 |\nabla G|\leq aG^{(n-1)/(n-p)},
 \qquad
 |\nabla\log G|\leq\frac{a}{b}.
\]
Both equalities hold for the normalized Euclidean
fundamental solution \(G(x)=|x-o|^{-a}\), so the constants
are sharp. The normalization matters: replacing \(G\)
by \(cG\) replaces \(b\) by \(c^{-1/a}b\) and changes
\(|\nabla b|\). The local estimate in
Lemma~\ref{lem:rough} supplies boundedness at infinity;
the normalized pole and the harmonic comparison function
\(G\) determine the constant \(1\). This is a global estimate for
the normalized minimal Green function. The local estimate of
\cite{WZ} applies to general positive \(p\)-harmonic functions;
here it supplies the preliminary bound needed for the global
comparison.

The comparison also has an equivalent conjugated form.
For a smooth function \(f\) on \(w>0\), the product rule
and \(\LG G=0\) give
\begin{equation}\label{eq:ground-state}
 G^{-2}\Div\bigl(G^2|\nabla G|^{p-2}\A\nabla f\bigr)
 =G^{-1}\LG(Gf).
\end{equation}
Indeed, both sides equal
\[
 \LG f+
 2G^{-1}|\nabla G|^{p-2}
       \langle\A\nabla G,\nabla f\rangle.
\]
Taking \(f=G^{-1}=b^a\) makes the right-hand side zero.
Thus \(b^a\) is an exact barrier for the conjugated operator.
The direct comparison above uses the same cancellation
through \(\LG G=0\).
\end{remark}

\subsection{Extension across the critical set}

The maximum-principle proof used the operator only on
compact regions where \(w>0\). Integration requires
control near \(Z\). We first establish the weighted
Hessian estimate that controls the error from a cutoff
around the critical set.

\begin{lemma}[A weighted Hessian estimate]
\label{lem:fixed-weighted}
Let \(u\) be \(p\)-harmonic on an open set \(U\), with
\(\Ric\geq0\). For every \(K\Subset U\),
\[
 \int_{K\cap\{|\nabla u|>0\}}
       |\nabla u|^{2p-4}|\Hess u|^2\dd\mu<\infty.
\]
\end{lemma}

\begin{proof}
Choose a smooth domain compactly contained in \(U\)
and containing \(K\). Solve the regularized Dirichlet problem
with boundary data \(u\):
\[
 \Div\bigl(t_\varepsilon^{(p-2)/2}
                 \nabla u_\varepsilon\bigr)=0,
 \qquad
 t_\varepsilon=\varepsilon^2+|\nabla u_\varepsilon|^2.
\]
The interior \(C^{1,\alpha}\) estimates for this regularization
\cite{Lewis,Tolksdorf} give, for every compact subset \(K'\) of
the chosen domain,
\[
 \|u_\varepsilon\|_{C^{1,\alpha}(K')}\leq C_{K'},
\]
where \(\alpha>0\) and \(C_{K'}\) are independent of
\(0<\varepsilon\leq1\); see also the regularization following problem~(A.1) in
\cite[Appendix~A]{BFMmass}. These are uniform estimates
for the regularized solutions, not merely regularity
of the limiting solution. The maximum principle bounds
\(\|u_\varepsilon\|_{L^\infty}\) by the fixed boundary
data. The structural ellipticity constants, after factoring
out \(t_\varepsilon^{(p-2)/2}\), lie between
\(\min\{1,p-1\}\) and \(\max\{1,p-1\}\), independently
of \(\varepsilon\). The constants in the interior estimate
may depend on the compact set, the fixed local geometry,
and the boundary data. The regularized
solutions minimize the strictly convex energies
\(p^{-1}\int(\varepsilon^2+|\nabla f|^2)^{p/2}\dd\mu\)
with boundary data \(u\). Variational compactness and uniqueness
of the original Dirichlet problem identify their weak
\(W^{1,p}\) limit with \(u\). The uniform H\"older bound on
the gradients then gives, after passing to a subsequence,
\(u_\varepsilon\to u\) in \(C^1_{\mathrm{loc}}\).
On compact subsets of \(\{|\nabla u|>0\}\), the equations
are uniformly elliptic for small \(\varepsilon\), and interior
estimates give smooth convergence.

For the following calculation, write \(t=t_\varepsilon\)
and define
\[
 \A_\varepsilon
 =g^{-1}+(p-2)
       \frac{\nabla u_\varepsilon\otimes\nabla u_\varepsilon}{t},
 \qquad
 L_\varepsilon f
 =\Div\bigl(t^{(p-2)/2}\A_\varepsilon\nabla f\bigr).
\]
Differentiating the regularized equation and commuting
covariant derivatives gives
\[
 \nabla_j\left(
       t^{(p-2)/2}\A_\varepsilon^{jk}
       \nabla_k\nabla_i u_\varepsilon\right)
 =t^{(p-2)/2}\Ric_i{}^\ell
       \nabla_\ell u_\varepsilon.
\]
Multiplying by \(\nabla^i u_\varepsilon\), summing over \(i\),
and using the product rule yields
\[
 \frac12L_\varepsilon t
 =t^{(p-2)/2}
 \left[
   |\Hess u_\varepsilon|^2
   +\frac{p-2}{t}
       |\Hess u_\varepsilon\nabla u_\varepsilon|^2
   +\Ric(\nabla u_\varepsilon,\nabla u_\varepsilon)
 \right].
\]

Apply the scalar chain rule to \(t^{p/2}\):
\begin{align*}
 \frac1pL_\varepsilon(t^{p/2})
 &=\frac12t^{(p-2)/2}L_\varepsilon t
   +\frac{p-2}{4}t^{p-3}
       |\nabla t|_{\A_\varepsilon}^2.
\end{align*}
Since
\[
 \nabla t=2\Hess u_\varepsilon\nabla u_\varepsilon,
\]
the Hessian terms combine as follows:
\begin{align*}
 |\Hess u_\varepsilon|_{\A_\varepsilon}^2
 ={}&|\Hess u_\varepsilon|^2
   +\frac{2(p-2)}{t}
       |\Hess u_\varepsilon\nabla u_\varepsilon|^2\\
 &+\frac{(p-2)^2}{t^2}
       \Hess u_\varepsilon
       (\nabla u_\varepsilon,\nabla u_\varepsilon)^2.
\end{align*}
We therefore obtain the exact identity
\[
 \frac1pL_\varepsilon(t^{p/2})
 =t^{p-2}\left[
       |\Hess u_\varepsilon|_{\A_\varepsilon}^2
       +\Ric(\nabla u_\varepsilon,\nabla u_\varepsilon)
       \right].
\]
In particular, the full Hessian expression is nonnegative
also when \(1<p<2\).

Every eigenvalue of \(\A_\varepsilon\) lies between
\(\min(1,p-1)\) and \(\max(1,p-1)\). Hence
\[
 |\Hess u_\varepsilon|_{\A_\varepsilon}^2
 \geq c(p)|\Hess u_\varepsilon|^2,
 \qquad c(p)=\min(1,p-1)^2>0.
\]
Choose a compactly supported cutoff \(\eta\) equal to \(1\)
on \(K\). Multiply the identity by \(\eta^2\), integrate
by parts, and use \(\Ric\geq0\).
Because
\[
 \nabla(t^{p/2})
 =pt^{(p-2)/2}\Hess u_\varepsilon\nabla u_\varepsilon,
\]
we obtain
\begin{align*}
 c(p)\int\eta^2t^{p-2}|\Hess u_\varepsilon|^2\dd\mu
 &\leq C(p)\int
       \eta|\nabla\eta|t^{p-3/2}
       |\Hess u_\varepsilon|\dd\mu\\
 &\leq\frac{c(p)}2
       \int\eta^2t^{p-2}|\Hess u_\varepsilon|^2\dd\mu+C(p)\int t^{p-1}|\nabla\eta|^2\dd\mu.
\end{align*}
Absorbing the first term on the right gives
\begin{equation}\label{eq:fixed-Caccioppoli}
 \int\eta^2t^{p-2}|\Hess u_\varepsilon|^2\dd\mu
 \leq C(p)\int t^{p-1}|\nabla\eta|^2\dd\mu.
\end{equation}
The right-hand side is uniformly bounded by the interior
gradient bounds and \(p-1>0\). Smooth convergence away
from the critical set and Fatou's lemma prove the assertion.
\end{proof}

Regularization is used here to obtain integrability.
The exact harmonic cancellation is retained for the
original solution. Indeed,
\[
 \A_\varepsilon\nabla u_\varepsilon
 =\left[(p-1)-\frac{(p-2)\varepsilon^2}{t}\right]
       \nabla u_\varepsilon,
\]
and therefore
\[
 L_\varepsilon u_\varepsilon
 =-(p-2)\varepsilon^2
       \Div\bigl(t^{(p-4)/2}\nabla u_\varepsilon\bigr).
\]
Thus the regularized linearization need not annihilate
\(u_\varepsilon\). We use the preceding estimate to pass
through the critical set of \(G\), while keeping the exact
identities for \(G\) itself.

\begin{proposition}[Weak defect identities]
\label{prop:fixed-weak}
Equations \eqref{eq:fixed-master},
\eqref{eq:fixed-div}, and \eqref{eq:second-basic}
hold in distributions on \(M\setminus\{o\}\), with no
additional measure supported on \(Z\).
The defect and the fluxes appearing in these identities
are locally integrable and are extended by zero on \(Z\).
At critical points, the distributional operators are understood
through these complete fluxes.
\end{proposition}

\begin{proof}
Fix a compact set away from \(o\), and write
\(s=|\nabla G|\). On a neighborhood of this set,
\(b\) and \(G\) are bounded above and away from zero.
The relation
\[
 w=\frac1a b^{a+1}s
\]
shows that \(w\) and \(s\) have the same zero set.
On \(s>0\), differentiation gives
\[
 |\nabla w|\leq C(|\Hess G|+s^2),
 \qquad
 |\Hess b|+|\Hess v|
 \leq C(|\Hess G|+s^2).
\]

Consider the flux
\[
 Y=b^{2-n-a}w^{p-2}\A\nabla F,
\]
extended by zero on \(Z\). Since
\(\nabla F=pw^{p-1}\nabla w\), the preceding bounds give
\[
 |Y|\leq C\bigl(s^{2p-3}|\Hess G|+s^{2p-1}\bigr).
\]
The definition of \(B\), the relation
\(\Delta_pv=nq^{p-1}F\), and boundedness of the Ricci tensor
on the compact set also give
\[
 \D\leq C\bigl(
       s^{2p-4}|\Hess G|^2+s^{2p}+s^{2p-2}\bigr).
\]
Lemma~\ref{lem:fixed-weighted}, applied to \(G\), controls
the first term in the last estimate. The remaining powers
are locally bounded because \(p>1\). Thus
\(\D\in L^1_{\mathrm{loc}}\).

For the flux estimate, write
\[
 s^{2p-3}|\Hess G|
 =\bigl(s^{p-2}|\Hess G|\bigr)s^{p-1}.
\]
The first factor is locally square integrable by the
weighted Hessian estimate, and the second is locally
bounded. Cauchy--Schwarz therefore gives
\(Y\in L^1_{\mathrm{loc}}\).

Choose a smooth function \(\chi\) satisfying
\[
 \chi=0\ \text{on }[0,1],\qquad
 \chi=1\ \text{on }[2,\infty),\qquad
 |\chi'|\leq C,
\]
and set \(\chi_\delta=\chi(s/\delta)\).
For \(\varphi\in C_c^\infty(M\setminus\{o\})\),
test the classical identity \eqref{eq:fixed-div}
on \(s>0\) with \(\varphi\chi_\delta\):
\begin{align*}
 -\int\chi_\delta
       \langle Y,\nabla\varphi\rangle\dd\mu
 -\int\varphi
       \langle Y,\nabla\chi_\delta\rangle\dd\mu
 =(p-1)q^{3-2p}
       \int\varphi\chi_\delta b^{-n-a}\D\dd\mu.
\end{align*}
Only the second term on the left needs further control.
It is supported in \(\{\delta<s<2\delta\}\), where
\[
 \delta^{-1}\leq\frac2s,
 \qquad |\nabla s|\leq|\Hess G|.
\]
Its absolute value is at most
\begin{align*}
 C\int_{\delta<s<2\delta}
       \left(s^{2p-4}|\Hess G|^2
             +s^{2p-2}|\Hess G|\right)\dd\mu\leq
 C'\int_{\delta<s<2\delta}
       \left(s^{2p-4}|\Hess G|^2+s^{2p}\right)\dd\mu.
\end{align*}
The last inequality follows from Young's inequality.
Its integrand is locally integrable, and the indicators
of these bands converge pointwise to zero.
Hence this error tends to zero.

Dominated convergence in the other terms gives
\[
 -\int\langle Y,\nabla\varphi\rangle\dd\mu
 =(p-1)q^{3-2p}
       \int\varphi b^{-n-a}\D\dd\mu.
\]
This proves \eqref{eq:fixed-div} in distributions and
excludes an additional measure on \(Z\).

Since \(b\in C^{1,\alpha}_{\mathrm{loc}}\) away from the pole,
\(b^{n+a-2}\) is an admissible \(C^1\) multiplier. Indeed,
the test function \(b^{n+a-2}\varphi\) can be approximated in
\(C^1\) by smooth compactly supported functions, and the
locally integrable flux and density permit passage to the limit.
The distributional product rule therefore gives
\begin{align*}
 \LB F
 =(p-1)q^{3-2p}b^{-2}\D
 +(p-1)(n+a-2)b^{-1}w^{p-2}
       \langle\nabla b,\nabla F\rangle,
\end{align*}
which is \eqref{eq:fixed-master}.

Finally, the flux of the product identity is, on \(s>0\),
\begin{align*}
 |\nabla G|^{p-2}\A\nabla(GF)
 =G|\nabla G|^{p-2}\A\nabla F+(p-1)F|\nabla G|^{p-2}\nabla G.
\end{align*}
It satisfies the same bound as \(Y\):
\[
 \left||\nabla G|^{p-2}\A\nabla(GF)\right|
 \leq C\bigl(s^{2p-3}|\Hess G|+s^{2p-1}\bigr).
\]
Testing \eqref{eq:second-basic} with
\(\varphi\chi_\delta\) therefore produces the same
vanishing cutoff error. Passing to the limit proves
the product identity in distributions.
\end{proof}

These weak identities can now be tested with functions
of \(b\) and integrated over annuli, as required for
the three formulas in Section~\ref{sec:fixed-formulas}.

\section{Three monotonicity formulas}
\label{sec:fixed-formulas}

Throughout this section, Assumption~\ref{ass:standing} is in
force. We retain the notation
\[
 a=\frac{n-p}{p-1},\qquad q=\frac{p}{p-1},
 \qquad w=|\nabla b|,\qquad F=w^p,
\]
and write
\[
 \Omega_r=\{b<r\},\qquad \Sigma_r=\{b=r\}.
\]
The quantities \(A_p\) and \(V_p\) are those defined in
\eqref{eq:fixed-AV}. We derive the three formulas separately.
Each uses a different form of the same Bochner identity:
\eqref{eq:fixed-master} gives the first formula,
the product identity \eqref{eq:second-basic} gives the second,
and the divergence identity \eqref{eq:fixed-div} gives the third.

The first step is to identify the appropriate level-set average.
For a regular level, the normalization of \(G\) gives
\[
 \int_{b=r}|\nabla G|^{p-1}\dd\sigma
 =a^{p-1}\om.
\]
Since
\[
 |\nabla G|^{p-1}
 =a^{p-1}r^{1-n}w^{p-1}
 \qquad\text{on }\Sigma_r,
\]
we obtain the conserved flux
\begin{equation}\label{eq:flux}
 r^{1-n}\int_{b=r}w^{p-1}\dd\sigma=\om.
\end{equation}
The same identity holds for almost every level by the weak
equation for \(G\).

This conservation law suggests defining
\begin{equation}\label{eq:I-def}
 I_u(r):=r^{1-n}\int_{b=r}u\,w^{p-1}\dd\sigma
 =a^{1-p}\int_{b=r}u|\nabla G|^{p-1}\dd\sigma.
\end{equation}
Thus \(I_1=\om\) and \(I_F=A_p\). In particular, the surface
exponent \(2p-1\) in \(A_p\) is determined by averaging the
Bochner quantity \(F=w^p\) against the conserved flux.
Up to the factor \(\om\), this is the surface average in
\cite[Lemma~5.3]{MRS}, specialized to the Euclidean model.
We repeat the differentiation calculation to make its relation
with \(\LG\) explicit. Initially, \(u\) is smooth on a band
of regular levels.

To compute its derivative, observe that
\[
 \Div\bigl(b^{1-n}w^{p-2}\nabla b\bigr)=0
 \qquad\text{away from }o.
\]
Indeed, the vector field inside the divergence equals
\(-a^{1-p}|\nabla G|^{p-2}\nabla G\).
On an annulus foliated by regular levels, the divergence theorem
and coarea give
\begin{align*}
 I_u(r)-I_u(s)
 =\int_{s<b<r}
       b^{1-n}w^{p-2}
       \langle\nabla u,\nabla b\rangle\dd\mu
 =\int_s^r t^{1-n}
       \int_{b=t}w^{p-2}\partial_\nu u
       \dd\sigma\dd t.
\end{align*}
Consequently,
\begin{equation}\label{eq:I-deriv}
 I_u'(r)
 =r^{1-n}\int_{b=r}w^{p-2}\partial_\nu u\dd\sigma.
\end{equation}

The same average is compatible with the linearization at \(G\).
Since
\[
 |\nabla G|^{p-2}
 =a^{p-2}r^{a+2-n}w^{p-2},
 \qquad \A\nu=(p-1)\nu
 \qquad\text{on }\Sigma_r,
\]
its boundary flux is
\begin{align*}
 \int_{b=r}|\nabla G|^{p-2}
       \langle\A\nabla u,\nu\rangle\dd\sigma
 =(p-1)a^{p-2}r^{a+2-n}
       \int_{b=r}w^{p-2}\partial_\nu u\dd\sigma
 =(p-1)a^{p-2}r^{a+1}I_u'(r).
\end{align*}
It follows that
\begin{equation}\label{eq:I-LG}
 [t^{a+1}I_u'(t)]_{t=s}^{t=r}
 =\frac{1}{(p-1)a^{p-2}}
       \int_{s<b<r}\LG u\dd\mu.
\end{equation}
The applications below use \(u=F\) and \(u=GF\).
For these functions, the estimates and weak identities established
in Section~\ref{sec:critical} justify the flux computations across
the critical set.

Before proving the formulas, we record the regularity and the
boundary limits needed in their proofs.

\begin{lemma}[Regularity and normalization]
\label{lem:fixed-averages}
The function \(A_p\) has a
\(W^{2,1}_{\mathrm{loc}}(0,\infty)\), hence \(C^1\),
representative. With this representative, \(V_p\) is \(C^1\) and
\begin{equation}\label{eq:fixed-V-deriv}
 V_p(r)=r^{-n}\int_0^r t^{n-1}A_p(t)\dd t,
 \qquad
 rV_p'(r)=A_p(r)-nV_p(r).
\end{equation}
Moreover,
\[
 0\leq A_p(r)\leq\om,\qquad
 0\leq V_p(r)\leq\frac{\om}{n},
\]
and
\[
 \lim_{r\downarrow0}A_p(r)=\om,\qquad
 \lim_{r\downarrow0}V_p(r)=\frac{\om}{n},\qquad
 \lim_{r\downarrow0}rV_p'(r)=0.
\]
There is a sequence \(s_j\downarrow0\) such that
\begin{equation}\label{eq:fixed-pole-sequence}
 s_jA_p'(s_j)\longrightarrow0.
\end{equation}
\end{lemma}

\begin{proof}
We first justify the absolute continuity of the flux
average across critical levels. Fix a compact band away
from the pole, put \(s=|\nabla G|\), and set
\[
 J:=b^{1-n}w^{p-2}\nabla b
   =-a^{1-p}|\nabla G|^{p-2}\nabla G.
\]
Thus \(\Div J=0\) weakly, and \(J\) extends continuously
by zero on \(Z\). Since \(w=a^{-1}b^{a+1}s\), on
\(s>0\) we have
\[
 |J|\,|\nabla F|
 \leq C\left(s^{2p-2}|\Hess G|+s^{2p}\right).
\]
The first term is integrable: write
\(s^{2p-2}|\Hess G|=(s^{p-2}|\Hess G|)s^p\)
and use Lemma~\ref{lem:fixed-weighted} and the local
gradient bound. The second term is locally bounded.
All densities involving derivatives are extended by zero
on \(Z\).

Let \(\chi_\delta=\chi(s/\delta)\) be the cutoff used in
Proposition~\ref{prop:fixed-weak}. Testing \(\Div J=0\)
with \(\varphi F\chi_\delta\), where \(\varphi\) is smooth
and compactly supported in the band, introduces the error
\[
 \left|\int \varphi F
          \langle J,\nabla\chi_\delta\rangle\dd\mu\right|
 \leq C\int_{\delta<s<2\delta}
                 s^{2p-2}|\Hess G|\dd\mu
 \longrightarrow0.
\]
The other terms converge by dominated convergence, giving
\(\Div(FJ)=\langle J,\nabla F\rangle\) in distributions.
Test this identity with \(\eta(b)\), where
\(\eta\in C_c^\infty(0,\infty)\). Properness of \(b\)
makes this test compactly supported away from \(o\),
and its \(C^1\) regularity suffices by smooth approximation.
Coarea yields
\[
 -\int_0^\infty A_p(r)\eta'(r)\dd r
 =\int_0^\infty\eta(r)r^{1-n}
       \int_{b=r}w^{p-2}\partial_\nu F\dd\sigma\dd r.
\]
The density on the right is locally integrable, since
its absolute integral is bounded by the corresponding
integral of \(|J|\,|\nabla F|\). Also \(A_p\) is locally
bounded by the conserved flux and local boundedness of
\(F\). Hence \(A_p\in W^{1,1}_{\mathrm{loc}}(0,\infty)\)
and \eqref{eq:I-deriv} holds weakly for \(u=F\).
Replacing \(F\) by \(GF\) only adds bounded \(C^1\)
factors on the band and locally bounded terms. The same
argument therefore justifies the average \(I_{GF}\).

The boundary flux of the vector field in
\eqref{eq:fixed-div} is
\begin{align*}
 \int_{b=r}b^{2-n-a}w^{p-2}
       \langle\A\nabla F,\nu\rangle\dd\sigma
 =(p-1)r^{2-n-a}
       \int_{b=r}w^{p-2}\partial_\nu F\dd\sigma
 =(p-1)r^{1-a}A_p'(r).
\end{align*}
Testing the weak divergence identity with functions of \(b\)
therefore gives
\begin{equation}\label{eq:fixed-density}
 \frac{d}{dr}\bigl(r^{1-a}A_p'(r)\bigr)
 =q^{3-2p}r^{-n-a}
       \int_{b=r,\,w>0}\frac{\D}{w}\dd\sigma
\end{equation}
in distributions. The right-hand side is locally integrable in
\(r\): by coarea, its integral over a compact interval is the
corresponding annular integral of \(b^{-n-a}\D\).
Thus \(r^{1-a}A_p'\in W^{1,1}_{\mathrm{loc}}\), and hence
\(A_p\in W^{2,1}_{\mathrm{loc}}\).
The annular flux identities consequently extend from
almost every pair of endpoints to every \(0<s<r\): their
boundary expressions are continuous, and their locally
integrable densities have continuous radial primitives.

Coarea applied to the definition of \(V_p\) gives
\[
 \int_{\Omega_r}w^{2p}\dd\mu
 =\int_0^r\int_{b=t}w^{2p-1}\dd\sigma\dd t
 =\int_0^r t^{n-1}A_p(t)\dd t.
\]
Differentiating proves \eqref{eq:fixed-V-deriv}.

The sharp gradient estimate \(w\leq1\), together with
\eqref{eq:flux}, gives
\[
 A_p(r)
 =r^{1-n}\int_{b=r}w^p w^{p-1}\dd\sigma
 \leq\om.
\]
The bound for \(V_p\) follows from its coarea representation.
The pole asymptotics \(w\to1\) and the conserved flux imply
\(A_p(r)\to\om\). The same representation then gives
\(V_p(r)\to\om/n\), and
\eqref{eq:fixed-V-deriv} gives \(rV_p'(r)\to0\).

Finally, apply the mean value theorem to \(A_p\) on
\([\varepsilon,2\varepsilon]\). At some
\(s\in(\varepsilon,2\varepsilon)\),
\[
 sA_p'(s)
 =\frac{s}{\varepsilon}
       \bigl(A_p(2\varepsilon)-A_p(\varepsilon)\bigr).
\]
The right-hand side tends to zero as
\(\varepsilon\downarrow0\). Choosing a sequence of disjoint
such intervals proves \eqref{eq:fixed-pole-sequence}.
\end{proof}

We use these representatives throughout the section.
The pole sequence allows us to remove the inner boundary
without assuming Hessian asymptotics at \(o\).

\begin{theorem}[The first monotonicity formula \eqref{eq:fixed-M1}]
\label{thm:fixed-first}
For every \(r>0\),
\[
 [A_p-q(n-1)V_p]'(r)
 =q^{3-2p}r^{-n-1}
       \int_{\Omega_r}\D\dd\mu.
\]
In particular, \(A_p-q(n-1)V_p\) is nondecreasing.
\end{theorem}

\begin{proof}
The volume quantity \(V_p\) has no radial weight in its
integrand. To obtain the corresponding defect integral,
multiply \eqref{eq:fixed-master} by \(b^2\):
\[
 b^2\LB F
 -(p-1)(n+a-2)bw^{p-2}
       \langle\nabla b,\nabla F\rangle
 =(p-1)q^{3-2p}\D.
\]
Integration by parts on \(\{s<b<r\}\) gives
\begin{align*}
 \int_{s<b<r}b^2\LB F\dd\mu
 =(p-1)[t^{n+1}A_p'(t)]_{t=s}^{t=r}
 -2(p-1)\int_{s<b<r}
       bw^{p-2}\langle\nabla b,\nabla F\rangle\dd\mu.
\end{align*}
Here the boundary term follows from
\eqref{eq:I-deriv} and \(\A\nu=(p-1)\nu\).

To evaluate the remaining integral, use
\(\Delta_pb=(n-1)b^{-1}w^p\):
\[
 \Div(bw^{p-2}\nabla b)
 =b\Delta_pb+w^p=nw^p.
\]
Applying the divergence theorem to
\(Fbw^{p-2}\nabla b\), we obtain
\begin{align*}
 \int_{s<b<r}
       bw^{p-2}\langle\nabla b,\nabla F\rangle\dd\mu
 &
 =[t^nA_p(t)]_{t=s}^{t=r}
       -n\int_{s<b<r}w^{2p}\dd\mu\\
 &
 =[t^n(A_p(t)-nV_p(t))]_{t=s}^{t=r}\\
 &
 =[t^{n+1}V_p'(t)]_{t=s}^{t=r}.
\end{align*}
Combining these calculations yields
\[
 [t^{n+1}(A_p'(t)-(n+a)V_p'(t))]_{t=s}^{t=r}
 =q^{3-2p}\int_{s<b<r}\D\dd\mu.
\]
Thus the coefficient \(n+a\) is determined by the integrations
by parts. Since
\[
 n+a=\frac{p(n-1)}{p-1}=q(n-1),
\]
it is exactly the coefficient in the statement.

Set \(s=s_j\), where \(s_j\) is given by
Lemma~\ref{lem:fixed-averages}. The inner boundary terms vanish:
\[
 s_j^{n+1}A_p'(s_j)
 =s_j^n\bigl(s_jA_p'(s_j)\bigr)\longrightarrow0,
\]
and
\[
 s_j^{n+1}V_p'(s_j)
 =s_j^n\bigl(s_jV_p'(s_j)\bigr)\longrightarrow0.
\]
Since \(\D\geq0\), monotone convergence gives
\[
 r^{n+1}\bigl(A_p'(r)-q(n-1)V_p'(r)\bigr)
 =q^{3-2p}\int_{\Omega_r}\D\dd\mu.
\]
This also proves that the defect integral is finite.
Dividing by \(r^{n+1}\) proves the formula.
\end{proof}

\begin{theorem}[The second monotonicity formula \eqref{eq:fixed-M2}]
\label{thm:fixed-second}
For every \(r>0\),
\[
 [r^{-a}(A_p(r)-\om)]'
 =q^{3-2p}r^{-a-1}
       \int_{\Omega_r}b^{-n}\D\dd\mu.
\]
Equivalently, for \(V_{p,\infty}\) defined in
\eqref{eq:fixed-Vinf},
\[
 [A_p-aV_{p,\infty}]'(r)
 =\frac{q^{3-2p}}{r}
       \int_{\Omega_r}b^{-n}\D\dd\mu.
\]
Both quantities are nondecreasing.
\end{theorem}

\begin{proof}
For this formula, use the product \(GF\) and the linearization
at \(G\). The product identity \eqref{eq:second-basic} is
compatible with the flux relation \eqref{eq:I-LG}.
Since \(G=r^{-a}\) on \(\Sigma_r\),
\[
 I_{GF}(r)=r^{-a}A_p(r),
\]
and therefore
\[
 r^{a+1}I_{GF}'(r)=rA_p'(r)-aA_p(r).
\]
Substituting \eqref{eq:second-basic} into
\eqref{eq:I-LG} gives
\begin{equation}\label{eq:fixed-second-shell}
 [tA_p'(t)-aA_p(t)]_{t=s}^{t=r}
 =q^{3-2p}
       \int_{s<b<r}b^{-n}\D\dd\mu.
\end{equation}

Along the pole sequence,
\[
 s_jA_p'(s_j)-aA_p(s_j)\longrightarrow-a\om.
\]
Monotone convergence consequently yields
\begin{equation}\label{eq:fixed-interior}
 rA_p'(r)-a(A_p(r)-\om)
 =q^{3-2p}\int_{\Omega_r}b^{-n}\D\dd\mu.
\end{equation}
In particular, the weighted defect is integrable at the pole.
Multiplying by \(r^{-a-1}\) gives the first identity.

For the alternative form, use
\[
 V_{p,\infty}'(r)=\frac{A_p(r)-\om}{r}
\]
and divide \eqref{eq:fixed-interior} by \(r\).
The two forms thus express the same defect identity.

The logarithmic primitive also has a volume interpretation.
For \(r>1\), coarea and \eqref{eq:flux} give
\[
 V_{p,\infty}(r)
 =\int_{1<b<r}
       b^{-n}(w^{2p}-w^p)\dd\mu.
\]
For \(r<1\), this is interpreted as a signed annular integral.
The fixed reference level makes the definition meaningful
without requiring a rate of convergence at the pole.
\end{proof}

\begin{theorem}[The third monotonicity formula \eqref{eq:fixed-M3}]
\label{thm:fixed-third}
For \(0<s<r\),
\[
 [t^{1-a}A_p'(t)]_{t=s}^{t=r}
 =q^{3-2p}
       \int_{s<b<r}b^{-n-a}\D\dd\mu.
\]
Hence \(r^{1-a}A_p'(r)\) is nondecreasing.
\end{theorem}

\begin{proof}
The divergence identity \eqref{eq:fixed-div} has boundary flux
\((p-1)r^{1-a}A_p'(r)\), as computed in the proof of
Lemma~\ref{lem:fixed-averages}. Integrating its radial form
\eqref{eq:fixed-density} gives
\begin{align*}
 [t^{1-a}A_p'(t)]_{t=s}^{t=r}
 &=q^{3-2p}\int_s^r t^{-n-a}
       \int_{b=t,\,w>0}
       \frac{\D}{w}\dd\sigma\dd t\\
 &=q^{3-2p}
       \int_{s<b<r}b^{-n-a}\D\dd\mu.
\end{align*}
The last equality is coarea. The weak identity across \(Z\)
ensures that no additional term occurs.
\end{proof}

The third formula also determines the sign of \(A_p'\).
Its boundary value at infinity follows from boundedness,
rather than from an assumed asymptotic expansion.

\begin{theorem}[Monotonicity of \(A_p\) and \(V_p\) \eqref{eq:fixed-exterior}]
\label{thm:fixed-decreasing}
The functions \(A_p\), \(V_p\), and \(V_{p,\infty}\) are
nonincreasing. Moreover,
\[
 \lim_{r\to\infty}r^{1-a}A_p'(r)=0,
\]
and, for every \(r>0\),
\[
 A_p'(r)
 =-q^{3-2p}r^{a-1}
       \int_{b>r}b^{-n-a}\D\dd\mu.
\]
The exterior defect integral is finite.
\end{theorem}

\begin{proof}
Suppose that \(A_p'(r_0)>0\). By the third formula,
\[
 A_p'(t)\geq r_0^{1-a}A_p'(r_0)t^{a-1}
 \qquad(t\geq r_0).
\]
Since \(a>0\), integration gives
\[
 A_p(R)-A_p(r_0)
 \geq r_0^{1-a}A_p'(r_0)
       \frac{R^a-r_0^a}{a}.
\]
This contradicts \(A_p\leq\om\). Thus \(A_p'\leq0\).

The function \(r^{1-a}A_p'(r)\) is now nondecreasing and
nonpositive, so it has a finite limit \(\ell\leq0\).
If \(\ell<0\), there is \(c>0\) such that
\[
 A_p'(t)\leq-ct^{a-1}
\]
for all sufficiently large \(t\). Integration would force
\(A_p(R)\to-\infty\), contrary to \(A_p\geq0\).
Therefore \(\ell=0\).

Letting the outer radius tend to infinity in the third formula
and using monotone convergence gives
\[
 -r^{1-a}A_p'(r)
 =q^{3-2p}\int_{b>r}b^{-n-a}\D\dd\mu.
\]
This proves the exterior representation and the finiteness
of its integral.

Since \(A_p\) is nonincreasing, its coarea average satisfies
\[
 V_p(r)-\frac{A_p(r)}{n}
 =r^{-n}\int_0^r
       t^{n-1}\bigl(A_p(t)-A_p(r)\bigr)\dd t
 \geq0.
\]
Equation~\eqref{eq:fixed-V-deriv} then gives \(V_p'\leq0\).
Finally,
\[
 V_{p,\infty}'(r)=\frac{A_p(r)-\om}{r}\leq0.
\]
\end{proof}

\begin{remark}
When \(p=2\), we have \(q=2\), \(a=n-2\), and
\[
 A_2(r)=r^{1-n}\int_{b=r}w^3\dd\sigma,\qquad
 V_2(r)=r^{-n}\int_{\Omega_r}w^4\dd\mu.
\]
Moreover,
\[
 q^{3-2p}=\frac12,\qquad
 q(n-1)=2(n-1),\qquad
 b^{-n-a}=b^{2-2n}.
\]
The three formulas therefore recover
\cite[Theorems~2.4, 2.8, and 2.20]{Colding}.
The logarithmic form and the exterior representation recover
his Theorem~2.9 and Corollary~2.21, respectively.

The volume normalization used here is the one in Colding's
original formulas. Section~\ref{sec:beta} instead uses the
weighted volume normalization of Colding--Minicozzi. At
\(\beta=p\), its surface quantity agrees with \(A_p\), while
its volume quantity is
\[
 r^{2-n}\int_{\Omega_r}b^{-2}w^{2p}\dd\mu.
\]
Keeping these two volume normalizations explicit makes the
relation between the two sections precise.
\end{remark}

The equality cases will be treated in
Section~\ref{sec:applications}.
The first two formulas measure the defect in a sublevel set
containing the normalized pole; equality therefore gives
global Euclidean rigidity.
The exterior formula measures the defect outside a level,
and equality gives a truncated metric cone.
If the derivative of \(r^{1-a}A_p'(r)\) vanishes at a regular
level, then \(A_p'\) also vanishes there, so the same exterior
cone conclusion follows.

\section{Rigidity and consequences of the three formulas}
\label{sec:applications}

The formulas in Section~\ref{sec:fixed-formulas} express
derivatives as integrals of the same nonnegative defect.
Their equality cases depend on the region measured by that
integral. Vanishing of the exterior defect gives a metric
cone outside a level. Vanishing of an interior defect reaches
the normalized pole and gives global Euclidean rigidity.
We first prove the exterior statement and then treat the
other equality cases.
For exterior Green levels, the cone conclusion agrees
with the capacitary rigidity in \cite{BFM}; the pointwise
gradient rigidity is the theorem of \cite{MRS} recovered
above. The proofs here make explicit how equality in each
of Colding's three formulas forces the vanishing of \(B\)
on the corresponding interior or exterior region.

\begin{proposition}[Exterior cone rigidity]
\label{prop:cone}
Let \(r_0>0\) be a regular level. If \(A_p'(r_0)=0\), then
every connected component of \(\{b>r_0\}\) is a truncated
metric cone, with
\[
 g=d\rho^2+\rho^2\gamma,\qquad
 b=\kappa\rho,\qquad
 \rho>\frac{r_0}{\kappa},
 \qquad 0<\kappa\leq1.
\]
Here \(\kappa\) is constant on the component, and \(\gamma\)
is a smooth metric on a compact connected
\((n-1)\)-dimensional cross-section.

Conversely, if every component of \(\{b>r_0\}\) has this
form with the stated Green radius, then \(A_p\) is constant
for \(r\geq r_0\).
\end{proposition}

\begin{proof}
By the exterior representation \eqref{eq:fixed-exterior},
\[
 0=-A_p'(r_0)
 =q^{3-2p}r_0^{a-1}
       \int_{b>r_0}b^{-n-a}\D\dd\mu.
\]
The integrand is nonnegative. Moreover, \(\D\) is smooth
on \(w>0\). Hence
\[
 \D=0
 \qquad\text{on }\{b>r_0\}\setminus Z.
\]
Since \(\Ric\geq0\), the definition of \(\D\) implies
\(B=0\) there.

We first show that an exterior component contains no
critical points. For every vector \(X\), differentiation
of \(w=|\nabla b|\) gives
\[
 \Hess b(\nu,X)=X(w).
\]
The last two terms in \eqref{eq:B-T} cancel when one
argument is \(\nu\), so
\[
 B(\nu,\cdot)=q^{p-1}bw^{p-2}\,dw.
\]
Thus \(B=0\) implies \(dw=0\) wherever \(w>0\).

Let \(E\) be a connected component of \(\{b>r_0\}\).
It contains a point with \(w>0\). Indeed, otherwise
\(\nabla b=0\) throughout \(E\), so \(b\) would be constant
there. The set \(E\) has nonempty boundary because \(M\)
is connected and \(o\notin E\), and continuity gives
\(b=r_0\) on that boundary. This contradicts \(b>r_0\)
in \(E\).

Choose a connected component \(U\) of \(E\cap\{w>0\}\).
The preceding calculation gives
\[
 w=\kappa>0\qquad\text{on }U.
\]
The set \(U\) is also closed in \(E\). To see this, let
\(x\in E\) be a limit point of \(U\). Continuity gives
\(w(x)=\kappa>0\). A sufficiently small connected
neighborhood of \(x\) lies in \(E\cap\{w>0\}\) and
intersects \(U\), so it belongs to the same component.
Consequently \(x\in U\).

Since \(U\) is nonempty, open, and closed in the connected
set \(E\), we have \(U=E\). Therefore
\[
 w=\kappa>0\qquad\text{throughout }E.
\]
The sharp gradient estimate gives \(\kappa\leq1\).

We now recover the metric. Since \(B=0\), equation
\eqref{eq:B-T} yields
\[
 b\Hess b=\kappa^2g-db\otimes db.
\]
Set \(\rho=b/\kappa\). Then
\[
 |\nabla\rho|=1,\qquad
 \Hess\rho=\frac1\rho
       \bigl(g-d\rho\otimes d\rho\bigr).
\]
Along the flow of \(\nabla\rho\),
\[
 \frac{d}{dt}\rho=|\nabla\rho|^2=1.
\]
Every closed band in \(E\) with
\[
 \frac{r_0}{\kappa}<c\leq\rho\leq d<\infty
\]
is compact by properness of \(b\). The flow can therefore
be continued between any two levels above \(r_0/\kappa\).
Each trajectory meets a fixed reference level exactly once,
giving a product representation
\[
 E\cong(r_0/\kappa,\infty)\times\Sigma,
\]
where \(\Sigma\) is compact and connected.

In these coordinates, the unit gradient is orthogonal to
the levels, so
\[
 g=d\rho^2+g_\rho.
\]
For tangent vectors transported by the flow,
\[
 \frac{d}{d\rho}g_\rho
 =(\mathcal L_{\nabla\rho}g)|_{T\Sigma}
 =2\Hess\rho|_{T\Sigma}
 =\frac2\rho g_\rho.
\]
If \(\rho_1\) is the reference level, integration gives
\[
 g_\rho=\frac{\rho^2}{\rho_1^2}g_{\rho_1}.
\]
Defining \(\gamma=\rho_1^{-2}g_{\rho_1}\), we obtain
\[
 g=d\rho^2+\rho^2\gamma.
\]

Conversely, on a cone with \(b=\kappa\rho\),
\[
 |\nabla b|=\kappa,\qquad
 b\Hess b=\kappa^2g-db\otimes db.
\]
Thus \(B=0\). The radial Ricci curvature of a metric cone
is zero, and \(\nabla v\) is radial, so \(\D=0\).
If this holds on every exterior component,
\eqref{eq:fixed-exterior} gives \(A_p'(r)=0\) for
\(r\geq r_0\).
\end{proof}

For a connected exterior, the conserved flux also determines
the relation between the cone metric and the constant gradient.
On \(b=r\), the induced metric is
\((r/\kappa)^2\gamma\), and hence
\[
 \om
 =r^{1-n}\kappa^{p-1}
       \left(\frac r\kappa\right)^{n-1}
       \Vol_\gamma(\Sigma)
 =\kappa^{p-n}\Vol_\gamma(\Sigma).
\]
Therefore
\[
 \kappa^{n-p}
 =\frac{\Vol_\gamma(\Sigma)}{\om},
 \qquad
 A_p(r)=\om\kappa^p
 \qquad(r\geq r_0).
\]
If the exterior has several components, the conserved flux
is the sum of their component fluxes.

The cross-section of an exterior cone need not be round.
If \(\kappa=1\) on any component, however,
Theorem~\ref{thm:sharp-main} gives global Euclidean rigidity.
Thus the sharp gradient estimate also distinguishes the
Euclidean equality case within the exterior cone conclusion.

\begin{corollary}[Equality in the third formula]
\label{cor:fixed-third-rigidity}
Let \(r_0>0\) be a regular level. If
\[
 \left.\frac{d}{dr}
       \bigl(r^{1-a}A_p'(r)\bigr)\right|_{r=r_0}=0,
\]
then \(A_p'(r_0)=0\), and every connected component of
\(\{b>r_0\}\) has the cone structure of
Proposition~\ref{prop:cone}.

More generally, if \(s<r\) are regular levels and
\[
 s^{1-a}A_p'(s)=r^{1-a}A_p'(r),
\]
then their common value is zero,
\(A_p'(t)=0\) for every \(t\geq s\), and every connected
component of \(\{b>s\}\) is a truncated metric cone.
\end{corollary}

\begin{proof}
Properness of \(b\) and regularity of \(r_0\) give a compact
band around \(\Sigma_{r_0}\) without critical points.
The level-set integrals are therefore smooth near \(r_0\).
Differentiating \eqref{eq:fixed-M3} and applying coarea yields
\[
 \left.\frac{d}{dr}
       \bigl(r^{1-a}A_p'(r)\bigr)\right|_{r=r_0}
 =q^{3-2p}r_0^{-n-a}
       \int_{b=r_0}\frac{\D}{w}\dd\sigma.
\]
If the derivative vanishes, the smooth nonnegative
integrand must vanish everywhere on \(\Sigma_{r_0}\).
Thus \(B=0\) on that level. In particular,
\[
 B_{\nu\nu}
 =q^{p-1}r_0w^{p-2}\partial_\nu w=0,
\]
so \(\partial_\nu F=0\). Equation~\eqref{eq:I-deriv}
now gives
\[
 A_p'(r_0)
 =r_0^{1-n}\int_{b=r_0}
       w^{p-2}\partial_\nu F\dd\sigma=0.
\]
Proposition~\ref{prop:cone} applies.

For the second statement, the third monotonicity formula
makes \(t^{1-a}A_p'(t)\) nondecreasing. Equality at \(s\)
and \(r\) makes it constant on \([s,r]\).
Since \(s\) is regular, its derivative at \(s\) vanishes,
and the first part gives \(A_p'(s)=0\).
The same function is nonpositive by
Section~\ref{sec:fixed-formulas}. Being nondecreasing
from the value zero, it remains zero for every \(t\geq s\).
The exterior cone conclusion follows again from
Proposition~\ref{prop:cone}.
\end{proof}

For the first two formulas, the vanishing defect is measured
inside a level. The pole normalization then fixes the
constant gradient at its Euclidean value. The same
conclusion follows from equality in the model bounds
for \(A_p\) and \(V_p\).

\begin{corollary}[Euclidean rigidity]
\label{cor:global-rigidity}
If, at some \(r_0>0\), either
\[
 [A_p-q(n-1)V_p]'(r_0)=0
 \qquad\text{or}\qquad
 [r^{-a}(A_p-\om)]'(r_0)=0,
\]
then \((M,g)\) is Euclidean and
\[
 G(x)=|x-o|^{-a}.
\]
The second condition can equivalently be replaced by
\[
 [A_p-aV_{p,\infty}]'(r_0)=0.
\]

The same conclusion follows if any one of
\[
 A_p(r_0)=\om,\qquad
 V_p(r_0)=\frac{\om}{n},\qquad
 V_p'(r_0)=0
\]
holds at a positive radius.

Conversely, all these equalities hold on Euclidean space
at every positive radius.
\end{corollary}

\begin{proof}
The first or second monotonicity formula gives
\(\D=0\) on the noncritical part of \(\Omega_{r_0}\).
On a sufficiently small connected punctured neighborhood
of \(o\), the pole asymptotics ensure \(w>0\).
There \(B=0\), and
\[
 B(\nu,\cdot)=q^{p-1}bw^{p-2}\,dw
\]
shows that \(w\) is constant. Since \(w\to1\) at the pole,
this constant is \(1\).
Theorem~\ref{thm:sharp-main} then gives the stated global
rigidity. The logarithmic version has the same weighted
defect integral, so its equality case is identical.

Suppose next that \(A_p(r_0)=\om\).
Since \(A_p\) is nonincreasing and bounded above by \(\om\),
\[
 A_p(t)=\om\qquad(0<t\leq r_0).
\]
Choose a sufficiently small regular level \(s<r_0\).
Subtracting the definition of \(A_p(s)\) from the conserved
flux gives
\[
 0=\om-A_p(s)
 =s^{1-n}\int_{b=s}
       w^{p-1}(1-w^p)\dd\sigma.
\]
On this level, \(w>0\) and \(w\leq1\). The integrand is
continuous and nonnegative, so \(w=1\) everywhere on
\(\Sigma_s\). The sharp gradient rigidity applies.

If \(V_p(r_0)=\om/n\), the coarea representation gives
\[
 0=\frac{\om}{n}-V_p(r_0)
 =r_0^{-n}\int_0^{r_0}
       t^{n-1}\bigl(\om-A_p(t)\bigr)\dd t.
\]
Nonnegativity and continuity imply \(A_p(t)=\om\)
for \(0<t\leq r_0\), reducing this case to the preceding one.

Finally, if \(V_p'(r_0)=0\), then
\eqref{eq:fixed-V-deriv} gives
\(V_p(r_0)=A_p(r_0)/n\). Consequently,
\[
 0=r_0^{-n}\int_0^{r_0}
       t^{n-1}\bigl(A_p(t)-A_p(r_0)\bigr)\dd t.
\]
The integrand is nonnegative because \(A_p\) is
nonincreasing. Hence \(A_p(t)=A_p(r_0)\) for
\(0<t\leq r_0\). Taking the pole limit gives
\(A_p(r_0)=\om\), and the previous argument applies.

On Euclidean space, \(b=|x-o|\) and \(w=1\), so
\[
 A_p=\om,\qquad V_p=\frac{\om}{n},\qquad
 V_{p,\infty}=0,\qquad \D=0.
\]
All the asserted equalities follow.
\end{proof}

The same sharp gradient bound and conserved flux also give
the capacity and geometric comparisons below.

\section{A family of monotonicity formulas}
\label{sec:beta}

We now extend the preceding formulas by varying the power
of the gradient, following the approach of
Colding--Minicozzi \cite{CM}. The Green radius, the linearized
operators, and the tensor \(B\) remain those introduced
above. The additional parameter changes the coefficients
of the normal and tangential components of \(B\).

Throughout this section, \(p\) is fixed. We use only
\(\beta\) as a subscript for the parameter family.
For \(\beta>0\), define
\begin{equation}\label{eq:AV-def}
 \begin{split}
 A_\beta(r)
 &=r^{1-n}\int_{b=r}w^{p-1+\beta}\dd\sigma,\\
 V_\beta(r)
 &=r^{2-n}\int_{b<r}b^{-2}w^{p+\beta}\dd\mu.
 \end{split}
\end{equation}
The surface quantity is the conserved flux average of
\(w^\beta\). The volume quantity has the weight used in
\cite{CM}.
At \(\beta=p\), \(A_\beta=A_p\); the weighted quantity
\(V_\beta\) differs from \(V_p\) in \eqref{eq:fixed-AV}.
Coarea gives
\begin{equation}\label{eq:V-deriv}
 V_\beta(r)
 =r^{2-n}\int_0^r t^{n-3}A_\beta(t)\dd t,
 \qquad
 rV_\beta'(r)=A_\beta(r)-(n-2)V_\beta(r).
\end{equation}
This derivative relation explains the normalization of
\(V_\beta\).

All level-set integrals are initially understood for almost
every level, with integration over its noncritical part.
The continuous representatives used below will be obtained
from the weak divergence identities.

\begin{lemma}\label{lem:beta-bounds}
Under Assumption~\ref{ass:standing}, for every \(\beta>0\),
the following bounds hold for almost every \(r>0\):
\[
 0\leq A_\beta(r)\leq\om,
 \qquad
 0\leq V_\beta(r)\leq\frac{\om}{n-2}.
\]
Moreover,
\[
 \lim_{r\downarrow0}A_\beta(r)=\om,
 \qquad
 \lim_{r\downarrow0}V_\beta(r)=\frac{\om}{n-2}.
\]
The surface limit is taken through radii where the
integral is defined. For \(\beta>\bs\), the bounds and limits
hold for the continuous representatives at all radii.
\end{lemma}

\begin{proof}
The sharp gradient estimate and the conserved flux give
\[
 A_\beta(r)
 =r^{1-n}\int_{b=r}w^\beta w^{p-1}\dd\sigma
 \leq r^{1-n}\int_{b=r}w^{p-1}\dd\sigma
 =\om.
\]
The bound for \(V_\beta\) follows from
\eqref{eq:V-deriv}. Near the pole, \(w\to1\) uniformly,
so the same flux identity gives \(A_\beta\to\om\).
Substitution into the coarea expression for \(V_\beta\)
gives its pole limit.
\end{proof}

\subsection{Differential identities}

The parameter identities follow from the scalar chain rule
applied to the fixed-exponent identity. Since
\(F=w^p\), we have
\[
 \nabla w^\beta
 =\frac{\beta}{p}w^{\beta-p}\nabla F.
\]
The chain rule for \(\LB\) therefore gives
\begin{align}
 \LB(w^\beta)
 &=
 \frac{\beta}{p}w^{\beta-p}\LB F
 +\frac{\beta(\beta-p)}{p^2}
       w^{\beta-2p}w^{p-2}|\nabla F|_{\A}^2
 \nonumber\\
 &=
 \frac{\beta}{p}w^{\beta-p}\LB F
 +\beta(\beta-p)w^{p+\beta-4}
                         |\nabla w|_{\A}^2.
 \label{eq:chain-beta}
\end{align}
The drift term in \eqref{eq:fixed-master} obeys the
first-order chain rule. Thus
\begin{align*}
 &\LB(w^\beta)
 -(p-1)(n+a-2)b^{-1}w^{p-2}
                  \langle\nabla b,\nabla w^\beta\rangle\\
 &\qquad
 =\beta q^{2-2p}b^{-2}w^{\beta-p}\D
  +\beta(\beta-p)w^{p+\beta-4}|\nabla w|_{\A}^2.
\end{align*}
For \(\beta>0\), on the regular set \(\{w>0\}\), define
\begin{equation}\label{eq:E-definition}
 \mathcal E_{p,\beta}
 :=q^{2-2p}b^{-2}w^{\beta-p}\mathcal D_p
   +(\beta-p)w^{p+\beta-4}|\nabla w|_{\A}^2,
\end{equation}
where \(w=|\nabla b|\), \(q=p/(p-1)\), and
\(\mathcal D_p\) is the defect defined in \eqref{Defect}.
The second term may be negative when \(\beta<p\).
To determine the correct parameter range, we express this
density through the components of \(B\).

The component calculation in Section~\ref{sec:operators}
gives the covector identity \eqref{Bnu}, i.e.
\[
 B(\nu,\cdot)=q^{p-1}bw^{p-2}\,dw.
\]
In particular,
\[
 q^{2-2p}b^{-2}w^{\beta-p}
                 |B(\nu,\cdot)|_{\A}^2
 =w^{p+\beta-4}|\nabla w|_{\A}^2.
\]
The trace-free condition
\[
 \sum_{i=2}^nB_{ii}=-(p-1)B_{\nu\nu}
\]
and the tangential trace decomposition give
\begin{align*}
 |B|_{\A}^2
 =&q^{2p-2}b^2w^{2p-2}|\mathring h|^2
     +\frac{n}{n-1}(p-1)^2B_{\nu\nu}^2+2(p-1)|B(\nu,\cdot)^\top|^2\\
 =&q^{2p-2}b^2w^{2p-2}|\mathring h|^2
     +\frac{n(p-1)}{n-1}|B(\nu,\cdot)|_{\A}^2+\frac{(p-1)(n-2)}{n-1}
                         |B(\nu,\cdot)^\top|^2.
\end{align*}
Here
\[
 |B(\nu,\cdot)|_{\A}^2
 =(p-1)B_{\nu\nu}^2+|B(\nu,\cdot)^\top|^2.
\]
Combining these identities with the curvature term in
\(\D\), we obtain
\begin{equation}\label{eq:E-B}
\begin{split}
 \E={}&w^{p+\beta-2}
       \left(|\mathring h|^2+\Ric(\nu,\nu)\right)\\
 &+\frac{q^{2-2p}}{n-1}b^{-2}w^{\beta-p}
 \left(
 \bigl[(n-1)\beta-(n-p)\bigr]
                         |B(\nu,\cdot)|_{\A}^2
 +(p-1)(n-2)|B(\nu,\cdot)^\top|^2
 \right).
\end{split}
\end{equation}
Consequently, the natural threshold is
\[
 \bs=\frac{n-p}{n-1}.
\]
Under \(\Ric\geq0\), the density is nonnegative whenever
\(\beta\geq\bs\).

For later geometric applications, we also record the
equivalent expression
\begin{equation}\label{eq:E-def}
 \E=w^{p+\beta-2}\left[
 |\mathring h|^2+\Ric(\nu,\nu)
 +(\beta+p-2)|\nabla^\Sigma\log w|^2
 +\frac{\beta-\bs}{p-1}
       \left(H-\frac{n-1}{b}w\right)^2
 \right].
\end{equation}
Indeed, \eqref{eq:normal-relation} gives
\[
 |\nabla\log w|_{\A}^2
 =|\nabla^\Sigma\log w|^2
 +\frac1{p-1}
       \left(H-\frac{n-1}{b}w\right)^2.
\]
Using \eqref{eq:D-geometric}, the two coefficients become
\[
 2(p-1)+\beta-p=\beta+p-2,
 \qquad
 \frac n{n-1}+\frac{\beta-p}{p-1}
 =\frac{\beta-\bs}{p-1}.
\]
At the fixed exponent,
\[
 \mathcal E_{p,p}=q^{2-2p}b^{-2}\D.
\]

\begin{proposition}\label{prop:master}
For every \(\beta>0\), the following identities hold
on \(w>0\):
\begin{align}
 \LB(w^\beta)
 &-(p-1)(n+a-2)b^{-1}w^{p-2}
             \langle\nabla b,\nabla w^\beta\rangle
 =\beta\E,
 \label{eq:master}\\
 \LG(Gw^\beta)
 &=\beta a^{p-2}b^{2-n}\E,
 \label{eq:master-product}\\
 \Div\bigl(b^{2-n-a}w^{p-2}\A\nabla w^\beta\bigr)
 &=\beta b^{2-n-a}\E.
 \label{eq:master-div}
\end{align}
\end{proposition}

\begin{proof}
The first identity is the chain-rule calculation above.
Multiplying it by \(b^{2-n-a}\) and moving the drift
inside the divergence gives \eqref{eq:master-div}.
For the product identity, observe that
\[
 G^2|\nabla G|^{p-2}
 =a^{p-2}b^{2-n-a}w^{p-2}.
\]
Since \(\LG G=0\), the product rule gives
\begin{align*}
 \LG(Gw^\beta)
 &=G^{-1}\Div
       \bigl(G^2|\nabla G|^{p-2}\A\nabla w^\beta\bigr)\\
 &=a^{p-2}b^a
       \Div\bigl(b^{2-n-a}w^{p-2}\A\nabla w^\beta\bigr)\\
 &=\beta a^{p-2}b^{2-n}\E.
\end{align*}
\end{proof}

We next justify integration across the critical set.
The additional gradient exponent requires a corresponding
weighted Hessian estimate. Its strict coercivity explains
why the exact global formulas below are stated for
\(\beta>\bs\).

\begin{lemma}[Weighted Hessian estimate]
\label{lem:weighted-H}
Let \(u\) be \(p\)-harmonic on an open set \(U\) in a
manifold with \(\Ric\geq0\). If \(\beta>\bs\) and
\(K\Subset U\), then
\[
 \int_{K\cap\{|\nabla u|>0\}}
       |\nabla u|^{p+\beta-4}|\Hess u|^2\dd\mu<\infty.
\]
\end{lemma}

\begin{proof}
As in Section~\ref{sec:critical}, solve a regularized
Dirichlet problem on a smooth intermediate domain:
\[
 \Div(t_\varepsilon^{(p-2)/2}\nabla u_\varepsilon)=0,
 \qquad
 t_\varepsilon=\varepsilon^2+|\nabla u_\varepsilon|^2.
\]
Set
\[
 \A_\varepsilon
 =I+(p-2)
       \frac{\nabla u_\varepsilon\otimes\nabla u_\varepsilon}
            {t_\varepsilon},
 \qquad
 L_\varepsilon f
 =\Div(t_\varepsilon^{(p-2)/2}
                         \A_\varepsilon\nabla f).
\]
The local compactness and regularity properties are those
used in Section~\ref{sec:critical}; see
\cite{Lewis,Tolksdorf} and the regularization following
problem~(A.1) in \cite[Appendix~A]{BFMmass}.

For this calculation, write \(t=t_\varepsilon\) and
\(H=\Hess u_\varepsilon\). Differentiating the equation gives
\[
 \frac12L_\varepsilon t
 =t^{(p-2)/2}
 \left[
 |H|^2+\frac{p-2}{t}|H\nabla u_\varepsilon|^2
 +\Ric(\nabla u_\varepsilon,\nabla u_\varepsilon)
 \right].
\]
Applying the scalar chain rule to \(t^{\beta/2}\), we obtain
\begin{equation}\label{eq:regularized-bochner}
 \frac1\beta L_\varepsilon(t^{\beta/2})
 =t^{(p+\beta-4)/2}
 \left[
 Q_\varepsilon+
 \Ric(\nabla u_\varepsilon,\nabla u_\varepsilon)
 \right],
\end{equation}
where
\begin{equation}\label{eq:Qeps}
 Q_\varepsilon
 =|H|^2+\frac{p+\beta-4}{t}|H\nabla u_\varepsilon|^2
 +\frac{(p-2)(\beta-2)}{t^2}
       H(\nabla u_\varepsilon,\nabla u_\varepsilon)^2.
\end{equation}

At a point with nonzero gradient, choose \(e_1\) in its
direction and put
\[
 z=\frac{|\nabla u_\varepsilon|^2}{t},
 \qquad
 \ell=1+(p-2)z.
\]
Then \(0\leq z\leq1\), and the equation implies
\[
 \sum_{i=2}^nH_{ii}=-\ell H_{11}.
\]
If \(\mathring H_\top\) denotes the trace-free part of
the tangential block, then
\[
 |H|^2
 =|\mathring H_\top|^2
   +\frac{\ell^2}{n-1}H_{11}^2
   +2\sum_{i=2}^nH_{1i}^2+H_{11}^2.
\]
Substituting into \eqref{eq:Qeps} gives
\begin{equation}\label{eq:coercivity}
 Q_\varepsilon=|\mathring H_\top|^2
 +[2+(p+\beta-4)z]\sum_{i=2}^nH_{1i}^2
+\ell\left[
 \frac n{n-1}
 +\left(\beta-2+\frac{p-2}{n-1}\right)z
 \right]H_{11}^2.
\end{equation}
The mixed coefficient has endpoint values \(2\) and
\(p+\beta-2\). The bracket in the normal coefficient has
endpoint values \(n/(n-1)\) and \(\beta-\bs\).
Both are uniformly positive when \(\beta>\bs\).
Since \(\ell\) is bounded above and away from zero,
\[
 Q_\varepsilon\geq c(n,p,\beta)|H|^2.
\]
At a zero-gradient point, \(Q_\varepsilon=|H|^2\)
directly.

Multiply \eqref{eq:regularized-bochner} by
\(\eta^2\), where \(\eta\) is compactly supported, and
integrate by parts. Since
\[
 \nabla(t^{\beta/2})
 =\beta t^{(\beta-2)/2}H\nabla u_\varepsilon,
\]
the boundary term is bounded by
\[
 C\int\eta|\nabla\eta|
             t^{(p+\beta-3)/2}|H|\dd\mu.
\]
Young's inequality absorbs half of the Hessian term
and gives
\begin{equation}\label{eq:Caccioppoli}
 \int\eta^2t^{(p+\beta-4)/2}|H|^2\dd\mu
 \leq C(n,p,\beta)
       \int t^{(p+\beta-2)/2}|\nabla\eta|^2\dd\mu.
\end{equation}
Here \(p+\beta-2>0\), since
\[
 p+\bs-2=\frac{(p-1)(n-2)}{n-1}>0.
\]
The uniform local gradient bounds make the right-hand
side bounded independently of \(\varepsilon\).
Smooth convergence away from the critical set and
Fatou's lemma prove the claim.
\end{proof}

\begin{proposition}[Weak identities and level-set regularity]
\label{prop:weak-div}
Let \(\beta>\bs\). Extend \(\E\) by zero on \(Z=\{w=0\}\).
Then the three identities in Proposition~\ref{prop:master}
hold distributionally away from the pole. Moreover,
\(A_\beta\) has a
\(W^{2,1}_{\mathrm{loc}}(0,\infty)\), hence \(C^1\),
representative.
\end{proposition}

\begin{proof}
Fix a compact set away from the pole and put
\(s=|\nabla G|\). On a slightly larger compact set,
\[
 w=a^{-1}b^{a+1}s.
\]
The factor \(a^{-1}b^{a+1}\) is positive and \(C^{1,\alpha}\),
with bounded first derivatives. It is smooth on \(s>0\).
For
\[
 Y_\beta=b^{2-n-a}w^{p-2}\A\nabla w^\beta,
\]
extended by zero on \(Z\), direct differentiation gives
\begin{align*}
 |Y_\beta|
 &\leq C\left(
 s^{p+\beta-3}|\Hess G|+s^{p+\beta-1}\right),\\
 |\E|
 &\leq C\left(
 s^{p+\beta-4}|\Hess G|^2
 +s^{p+\beta-2}+s^{p+\beta}\right).
\end{align*}
To check the flux integrability, write
\[
 s^{p+\beta-3}|\Hess G|
 =\left(s^{(p+\beta-4)/2}|\Hess G|\right)
        s^{(p+\beta-2)/2}.
\]
The first factor lies in \(L^2_{\mathrm{loc}}\) by
Lemma~\ref{lem:weighted-H}; the second is locally bounded
because \(p+\beta-2>0\). Cauchy--Schwarz gives
\(Y_\beta\in L^1_{\mathrm{loc}}\), and the same weighted
estimate gives \(\E\in L^1_{\mathrm{loc}}\).

Use a cutoff \(\chi_\delta=\chi(s/\delta)\), equal
to zero for \(s\leq\delta\) and one for \(s\geq2\delta\).
Testing \eqref{eq:master-div} with
\(\varphi\chi_\delta\), the additional term is bounded by
\[
 C\int_{\delta<s<2\delta}
 \left(s^{p+\beta-4}|\Hess G|^2+s^{p+\beta}\right)\dd\mu.
\]
It tends to zero by the weighted estimate.
Thus \eqref{eq:master-div} extends across \(Z\)
without an additional measure there. Multiplication by the
locally \(C^1\) function \(b^{n+a-2}\), justified by smooth
\(C^1\) approximation of the test functions as in
Proposition~\ref{prop:fixed-weak}, gives \eqref{eq:master}.
For the product identity, its flux on \(s>0\) is
\[
 |\nabla G|^{p-2}\A\nabla(Gw^\beta)
 =G|\nabla G|^{p-2}\A\nabla w^\beta
   +(p-1)w^\beta|\nabla G|^{p-2}\nabla G.
\]
Its norm is bounded by
\(C(s^{p+\beta-3}|\Hess G|+s^{p+\beta-1})\), just as for
\(Y_\beta\). Testing \eqref{eq:master-product} with
\(\varphi\chi_\delta\) therefore gives the same vanishing
cutoff error. This proves its distributional extension.

To obtain the regularity of the average, use
\(J=b^{1-n}w^{p-2}\nabla b\) as in the proof of
Lemma~\ref{lem:fixed-averages}. On the compact set,
\[
 |J|\,|\nabla w^\beta|
 \leq C\left(s^{p+\beta-2}|\Hess G|
                    +s^{p+\beta}\right)\in L^1.
\]
Indeed,
\[
 s^{p+\beta-2}|\Hess G|
 =\left(s^{(p+\beta-4)/2}|\Hess G|\right)
        s^{(p+\beta)/2},
\]
and Lemma~\ref{lem:weighted-H}, local boundedness of
\(s\), and Cauchy--Schwarz apply. Testing \(\Div J=0\)
with \(\varphi w^\beta\chi_\delta\) produces an error
bounded by
\[
 C\int_{\delta<s<2\delta}
               s^{p+\beta-2}|\Hess G|\dd\mu
 \longrightarrow0.
\]
Thus \(\Div(w^\beta J)=\langle J,\nabla w^\beta\rangle\)
holds distributionally. Testing with \(\eta(b)\) and
using coarea, exactly as in Lemma~\ref{lem:fixed-averages},
gives a locally integrable weak derivative of the locally
bounded average. Therefore
\(A_\beta\in W^{1,1}_{\mathrm{loc}}\), with
\[
 A_\beta'(r)
 =r^{1-n}\int_{b=r}
                  w^{p-2}\partial_\nu w^\beta\dd\sigma
 \qquad\text{for almost every }r.
\]
Testing \eqref{eq:master-div} with functions of \(b\)
and applying coarea then gives
\[
 \frac d{dr}\bigl(r^{1-a}A_\beta'(r)\bigr)
 =\frac{\beta}{p-1}r^{2-n-a}
       \int_{b=r}\frac{\E}{w}\dd\sigma
\]
in distributions. Its right-hand side is locally
integrable in \(r\), because its integral on an interval
is an annular defect integral. Hence
\(r^{1-a}A_\beta'\in W^{1,1}_{\mathrm{loc}}\), which
implies the asserted regularity.

This argument uses coarea on the noncritical part of
the levels and does not require Sard's theorem for \(b\).
\end{proof}

For later use, on a regular level, \eqref{eq:I-deriv} and
\eqref{eq:normal-relation} give
\begin{equation}\label{eq:A-H}
\begin{split}
 A_\beta'(r)
 &=\beta r^{1-n}\int_{b=r}
       w^{p+\beta-3}\partial_\nu w\dd\sigma\\
 &=-\frac{\beta}{p-1}r^{1-n}
       \int_{b=r}w^{p+\beta-2}
          \left(H-\frac{n-1}{r}w\right)\dd\sigma.
\end{split}
\end{equation}

\subsection{The three monotonicity formulas}

We now prove the three formulas separately. Throughout
this subsection, \(\beta>\bs\), and the representatives
provided by Proposition~\ref{prop:weak-div} are understood.

\begin{theorem}[The first monotonicity formula]
\label{thm:beta-first}
For \(r>0\),
\begin{equation}\label{eq:M1}
 \bigl[A_\beta-(n+a-2)V_\beta\bigr]'(r)
 =\frac{\beta}{p-1}r^{1-n}
             \int_{b<r}\E\dd\mu.
\end{equation}
In particular, \(A_\beta-(n+a-2)V_\beta\) is
nondecreasing.
\end{theorem}

\begin{proof}
The boundary flux of \(\LB(w^\beta)\) on \(\Sigma_r\)
is
\begin{align*}
 \int_{b=r}w^{p-2}
                 \langle\A\nabla w^\beta,\nu\rangle\dd\sigma
 &=(p-1)\int_{b=r}
                 w^{p-2}\partial_\nu w^\beta\dd\sigma=(p-1)r^{n-1}A_\beta'(r).
\end{align*}
The drift term is related to \(V_\beta'\).
Indeed,
\[
 \Div(b^{-1}w^{p-2}\nabla b)
 =(n-2)b^{-2}w^p.
\]
Applying the product rule and integrating on
\(\{s<b<r\}\), we obtain
\begin{align*}
 \int_{s<b<r}b^{-1}w^{p-2}
                 \langle\nabla b,\nabla w^\beta\rangle\dd\mu
 =&[t^{n-2}A_\beta(t)]_s^r
 -(n-2)\int_{s<b<r}b^{-2}w^{p+\beta}\dd\mu\\
 =&[t^{n-1}V_\beta'(t)]_s^r.
\end{align*}
Integration of \eqref{eq:master} therefore gives
\[
 (p-1)
 \left[t^{n-1}
       \bigl(A_\beta'-(n+a-2)V_\beta'\bigr)\right]_s^r
 =\beta\int_{s<b<r}\E\dd\mu.
\]
Since \(A_\beta(r)\to\om\), there is a sequence
\(s_j\downarrow0\) such that
\(s_jA_\beta'(s_j)\to0\). Also,
\[
 rV_\beta'(r)
 =A_\beta(r)-(n-2)V_\beta(r)\longrightarrow0.
\]
Because \(n>2\), both inner boundary terms vanish
along this sequence. Letting \(s=s_j\downarrow0\)
proves \eqref{eq:M1}. Nonnegativity also proves
finiteness of its defect integral.
\end{proof}

By coarea, the derivative can equivalently be written as
\[
 \bigl[A_\beta-(n+a-2)V_\beta\bigr]'(r)
 =\frac{\beta}{p-1}r^{1-n}
       \int_0^r\int_{b=t}\frac{\E}{w}
                                      \dd\sigma\dd t.
\]
Together with \eqref{eq:E-B}, this is the direct
counterpart of the first formula in \cite{CM}.

\begin{theorem}[The second monotonicity formula]
\label{thm:beta-second}
For \(0<s<r\),
\begin{equation}\label{eq:M2}
 [tA_\beta'(t)-aA_\beta(t)]_s^r
 =\frac{\beta}{p-1}
       \int_{s<b<r}b^{2-n}\E\dd\mu.
\end{equation}
Moreover,
\begin{equation}\label{eq:interior}
 rA_\beta'(r)-a(A_\beta(r)-\om)
 =\frac{\beta}{p-1}
       \int_{b<r}b^{2-n}\E\dd\mu,
\end{equation}
and hence
\begin{equation}\label{eq:second-scaled}
 \bigl[r^{-a}(A_\beta(r)-\om)\bigr]'
 =\frac{\beta}{p-1}r^{-a-1}
       \int_{b<r}b^{2-n}\E\dd\mu.
\end{equation}
Thus \(r^{-a}(A_\beta-\om)\) is nondecreasing.
\end{theorem}

\begin{proof}
On a level of \(b\),
\[
 I_{Gw^\beta}(r)=r^{-a}A_\beta(r),
 \qquad
 r^{a+1}I_{Gw^\beta}'(r)
 =rA_\beta'(r)-aA_\beta(r).
\]
Using \eqref{eq:I-LG} and \eqref{eq:master-product},
we obtain
\begin{align*}
 [tA_\beta'(t)-aA_\beta(t)]_s^r
 =&\frac1{(p-1)a^{p-2}}
       \int_{s<b<r}\LG(Gw^\beta)\dd\mu\\
 =&\frac{\beta}{p-1}
       \int_{s<b<r}b^{2-n}\E\dd\mu.
\end{align*}
This proves \eqref{eq:M2}.

Along the pole sequence used above,
\[
 s_jA_\beta'(s_j)-aA_\beta(s_j)\longrightarrow-a\om.
\]
Letting \(s=s_j\downarrow0\) gives
\eqref{eq:interior}. Multiplication by \(r^{-a-1}\)
gives \eqref{eq:second-scaled}.
\end{proof}

\begin{theorem}[The third monotonicity formula]
\label{thm:beta-third}
For \(0<s<r\),
\begin{equation}\label{eq:M3}
 [t^{1-a}A_\beta'(t)]_s^r
 =\frac{\beta}{p-1}
       \int_{s<b<r}b^{2-n-a}\E\dd\mu.
\end{equation}
In particular, \(r^{1-a}A_\beta'(r)\) is
nondecreasing.
\end{theorem}

\begin{proof}
The boundary flux in \eqref{eq:master-div} is
\begin{align*}
 \int_{b=r}b^{2-n-a}w^{p-2}
                 \langle\A\nabla w^\beta,\nu\rangle\dd\sigma
 =(p-1)r^{2-n-a}
       \int_{b=r}w^{p-2}\partial_\nu w^\beta\dd\sigma
 =(p-1)r^{1-a}A_\beta'(r).
\end{align*}
Integration over \(\{s<b<r\}\), followed by division
by \(p-1\), proves the formula.
\end{proof}

The third formula controls the derivative of \(A_\beta\).
Its sign follows by combining this formula with the
uniform bound in Lemma~\ref{lem:beta-bounds}, as in
\cite{CM}.

\begin{theorem}[Monotonicity of \(A_\beta\)]
\label{thm:main}
The functions \(A_\beta\) and \(V_\beta\) are
nonincreasing. Moreover,
\begin{equation}\label{eq:exterior}
 A_\beta'(r)
 =-\frac{\beta}{p-1}r^{a-1}
       \int_{b>r}b^{2-n-a}\E\dd\mu.
\end{equation}
The exterior integral is finite for every \(r>0\).
\end{theorem}

\begin{proof}
The third formula makes \(r^{1-a}A_\beta'(r)\) nondecreasing.
If \(A_\beta'(r_0)>0\), then
\[
 A_\beta(R)-A_\beta(r_0)
 \geq r_0^{1-a}A_\beta'(r_0)\frac{R^a-r_0^a}{a},
 \qquad R\geq r_0,
\]
contradicting \(A_\beta\leq\om\). Thus \(A_\beta'\leq0\),
and the scaled derivative has a nonpositive limit at infinity.
A negative limit would force \(A_\beta\to-\infty\) by
integration, so the limit is zero. Letting the outer radius
tend to infinity in \eqref{eq:M3} proves \eqref{eq:exterior}
and finiteness of its integral.

Finally, the coarea average \eqref{eq:V-deriv} gives
\(V_\beta\geq A_\beta/(n-2)\), because \(A_\beta\) is
nonincreasing. The same identity then gives \(V_\beta'\leq0\).
\end{proof}

\begin{remark}[The harmonic case]
For \(p=2\),
\[
 a=n-2,\qquad n+a-2=2(n-2),
 \qquad \bs=\frac{n-2}{n-1}.
\]
The tensor \(\A\) is the identity, and \eqref{eq:E-B}
becomes
\begin{align*}
 \mathcal E_{2,\beta}
 =w^\beta\bigl(|\mathring h|^2+\Ric(\nu,\nu)\bigr)
 +\frac{b^{-2}w^{\beta-2}}{4(n-1)}
 \bigg(
 \bigl[1+(\beta-1)(n-1)\bigr]|B(\nu,\cdot)|^2
 +(n-2)|B(\nu,\cdot)^\top|^2
 \bigg).
\end{align*}
Thus the three formulas reduce respectively to
\cite[Theorems 3.9, 3.17, and 3.27]{CM}.
The exterior representation reduces to
\cite[Theorem 1.3]{CM}.
\end{remark}

The volume normalization in this section differs from
the unweighted \(V_p\) used earlier. To recover that
normalization, multiply \eqref{eq:master} by \(b^2\)
and repeat the integration in Section~\ref{sec:fixed-formulas}.
This gives
\[
 \left[
 A_\beta(r)
 -(n+a)r^{-n}\int_{b<r}w^{p+\beta}\dd\mu
 \right]'
 =\frac{\beta}{p-1}r^{-n-1}
       \int_{b<r}b^2\E\dd\mu.
\]
At \(\beta=p\), this is precisely
\eqref{eq:fixed-M1}. The second and third formulas
directly recover \eqref{eq:fixed-M2} and
\eqref{eq:fixed-M3}, using
\(\mathcal E_{p,p}=q^{2-2p}b^{-2}\D\).

At the endpoint \(\beta=\bs\), the pointwise identities
remain valid on \(w>0\), and \(\E\geq0\). However,
the normal coefficient in \eqref{eq:coercivity} loses
strict positivity. The preceding argument therefore
does not prove the exact global defect identities
at this exponent.

There are nevertheless endpoint monotonicity statements.
Let \(\beta_j\downarrow\bs\). Dominated convergence
with respect to the finite flux measure gives
\(A_{\beta_j}\to A_{\bs}\) almost everywhere and locally
in \(L^1\); the same holds for \(V_{\beta_j}\).
Consequently, \(A_{\bs}\) and \(V_{\bs}\) have
nonincreasing representatives, whereas
\[
 A_{\bs}-(n+a-2)V_{\bs},
 \qquad
 r^{-a}(A_{\bs}-\om)
\]
have nondecreasing representatives.
The third formula also gives convexity after the change
of variable \(t=r^a\), since
\[
 \frac{d}{dt}A_{\beta_j}(t^{1/a})
 =\frac1a r^{1-a}A_{\beta_j}'(r).
\]
Thus \(t\mapsto A_{\beta_j}(t^{1/a})\) is convex.
Its limit has a convex representative, so
\(r^{1-a}A_{\bs}'\) is nondecreasing distributionally.
These limiting statements do not establish an exact
absolutely continuous defect representation at the
endpoint. In the harmonic case, the endpoint formulas
are included in \cite{CM}.

\subsection{Rigidity}\label{sec:rigidity}
\begin{theorem}[Rigidity]
\label{thm:beta-rigidity}
Under Assumption~\ref{ass:standing}, let \(\beta>\bs\).

\begin{enumerate}
\item
If, at some \(r_0>0\), either
\[
 \bigl[A_\beta-(n+a-2)V_\beta\bigr]'(r_0)=0
\]
or
\[
 \left[r^{-a}(A_\beta(r)-\om)\right]'_{r=r_0}=0,
\]
then \(M\) is isometric to Euclidean space and
\(G(x)=|x-o|^{-a}\).

\item
If \(r_0>0\) is a regular value of \(b\) and
\(A_\beta'(r_0)=0\), then every connected component
of \(\{b>r_0\}\) is a truncated metric cone. More
precisely, on each component there are a constant
\(0<\kappa\leq1\) and a fixed metric \(\gamma\) on
a compact cross-section such that
\[
 g=d\rho^2+\rho^2\gamma,
 \qquad b=\kappa\rho,
 \qquad \rho>\frac{r_0}{\kappa}.
\]

\item
If \(r_0>0\) is a regular value and
\[
 \left.
 \frac d{dr}\left(r^{1-a}A_\beta'(r)\right)
 \right|_{r=r_0}=0,
\]
then the same exterior cone conclusion holds.

In particular, if \(0<s<r\) are regular values and
\[
 s^{1-a}A_\beta'(s)=r^{1-a}A_\beta'(r),
\]
then this common value is zero,
\(A_\beta'(t)=0\) for every \(t\geq s\), and every
connected component of \(\{b>s\}\) is a truncated
metric cone.
\end{enumerate}

Conversely, on a cone exterior with \(b=\kappa\rho\),
\(A_\beta\) is constant on the exterior scales.
On Euclidean space, all three monotone quantities
are constant.
\end{theorem}

\begin{proof}
The strict parameter condition gives
\[
 \beta-\bs>0,
 \qquad
 \beta+p-2>0.
\]
Thus the geometric expression \eqref{eq:E-def}
shows that, wherever \(w>0\), vanishing of \(\E\)
implies
\begin{equation}\label{eq:beta-equality-equations}
 \mathring h=0,\qquad
 \nabla^\Sigma w=0,\qquad
 H=\frac{n-1}{b}w,\qquad
 \Ric(\nu,\nu)=0.
\end{equation}
The normal relation \eqref{eq:normal-relation}
then gives
\[
 (p-1)\partial_\nu w
 =-w\left(H-\frac{n-1}{b}w\right)=0.
\]
Consequently, \(\nabla w=0\) on every connected
noncritical region where the defect vanishes.

For the first assertion, either
\eqref{eq:M1} or \eqref{eq:second-scaled} gives
\[
 \E=0
 \quad\text{on the noncritical part of }\{b<r_0\}.
\]
Here the conclusion is pointwise on the regular
set, since the density is smooth there and
nonnegative.

The pole asymptotics ensure that \(w>0\) in a
punctured neighborhood of \(o\). On a sufficiently
small connected punctured neighborhood, the
preceding equations make \(w\) constant. Its
pole limit is one, so \(w=1\) there.
The equality case of
Theorem~\ref{thm:sharp-main} therefore gives
global Euclidean rigidity and the stated
Green function.

For the second assertion, the exterior formula
\eqref{eq:exterior} gives
\[
 \E=0
 \qquad\text{on }\{b>r_0\}\setminus Z.
\]
Equations \eqref{eq:beta-equality-equations} and
\eqref{eq:D-geometric} then imply that \(\D=0\)
on this set. Hence \eqref{eq:fixed-exterior} yields
\(A_p'(r_0)=0\). Proposition~\ref{prop:cone}
gives the stated cone structure on each connected
component of \(\{b>r_0\}\).

For the third assertion, regularity of \(r_0\)
and properness of \(b\) give a smooth foliation
by nearby compact levels. The third formula
and coarea therefore yield the pointwise identity
\[
 \left.
 \frac d{dr}\left(r^{1-a}A_\beta'(r)\right)
 \right|_{r=r_0}
 =
 \frac{\beta}{p-1}r_0^{2-n-a}
       \int_{b=r_0}\frac{\E}{w}\dd\sigma.
\]
If the left-hand side is zero, nonnegativity
and \(w>0\) on \(\Sigma_{r_0}\) imply
\(\E=0\) there. Since \(\beta>\bs\), we obtain
\[
 H=\frac{n-1}{r_0}w
 \quad\text{on }\Sigma_{r_0}.
\]
The first-variation formula \eqref{eq:A-H}
then gives \(A_\beta'(r_0)=0\). The second
assertion applies.

Finally, suppose the scaled derivative has
equal values at the regular radii \(s<r\).
By its monotonicity, it is constant on \([s,r]\).
Its derivative at \(s\) is therefore zero,
and the third assertion gives
\(A_\beta'(s)=0\). Hence
\[
 s^{1-a}A_\beta'(s)=0.
\]
The scaled derivative is nondecreasing and
nonpositive by Theorem~\ref{thm:main}; it must
therefore vanish for every \(t\geq s\).
The exterior cone conclusion follows.

Conversely, on a cone with \(b=\kappa\rho\),
the gradient is constant, the levels are
umbilic, their mean curvature is
\((n-1)w/b\), and the radial Ricci curvature
vanishes. Thus \(\E=0\) there, and
\eqref{eq:A-H} gives \(A_\beta'=0\).
The Euclidean model has zero defect everywhere,
so all three monotone quantities are constant.
\end{proof}

\begin{corollary}[Rigidity in the Euclidean bounds]
\label{cor:beta-model-rigidity}
Under the same assumptions, if at some \(r_0>0\)
\[
 A_\beta(r_0)=\om,
 \qquad\text{or}\qquad
 V_\beta(r_0)=\frac{\om}{n-2},
 \qquad\text{or}\qquad
 V_\beta'(r_0)=0,
\]
then \(M\) is Euclidean.
\end{corollary}

\begin{proof}
If \(A_\beta(r_0)=\om\), bounded monotonicity gives
\(A_\beta(t)=\om\) for \(0<t\leq r_0\).

If \(V_\beta(r_0)=\om/(n-2)\), its coarea
representation and \(A_\beta\leq\om\) give the
same conclusion. If \(V_\beta'(r_0)=0\), then
\[
 0
 =V_\beta(r_0)-\frac{A_\beta(r_0)}{n-2}
 =r_0^{2-n}\int_0^{r_0}
       t^{n-3}[A_\beta(t)-A_\beta(r_0)]\dd t.
\]
The integrand is nonnegative. Continuity and
the pole limit imply
\(A_\beta(t)=\om\) for \(0<t\leq r_0\).

On a sufficiently small regular level,
the flux identity now gives
\[
 \int_{b=t}w^{p-1}(1-w^\beta)\dd\sigma=0.
\]
Since \(0<w\leq1\) there, we obtain \(w=1\).
Theorem~\ref{thm:sharp-main} proves the claim.
\end{proof}

\section{Geometric consequences and related monotonicity formulas}
\label{sec:comparison}

We first derive geometric estimates from the parameter family,
then record its relation to capacitary and scalar curvature
monotonicity. Unless stated otherwise,
Assumption~\ref{ass:standing} and \(\beta>\bs\) remain in force.
The representatives are those of Section~\ref{sec:beta};
integrals on nonsmooth levels are taken over their noncritical
parts and understood for almost every radius.

\subsection{Geometric consequences}

On a regular level, the first variation \eqref{eq:A-H} and
\(A_\beta'\leq0\) give the weighted mean-curvature comparison
\begin{equation}\label{eq:weighted-Minkowski}
 \int_{b=r}H w^{p+\beta-2}\dd\sigma
 \geq\frac{n-1}{r}
       \int_{b=r}w^{p+\beta-1}\dd\sigma.
\end{equation}
This compares the integrated mean curvature with the
cone model and requires no pointwise mean convexity.

For a bounded domain \(\Omega\) with smooth boundary,
we use the unnormalized \(p\)-capacity
\begin{equation}\label{eq:p-capacity}
 \Cap_p(\Omega)
 :=\inf\left\{
   \int_M|\nabla f|^p\dd\mu:
   f\in C_c^\infty(M),\
   f\geq1\ \text{on }\overline\Omega
 \right\}.
\end{equation}
If \(u\) is its exterior \(p\)-capacitary potential,
then
\[
 \Cap_p(\Omega)
 =\int_{M\setminus\overline\Omega}|\nabla u|^p\dd\mu
 =\int_{\partial\Omega}|\nabla u|^{p-1}\dd\sigma.
\]
For a regular Green level \(\Omega_r=\{b<r\}\),
the potential is \(u_r=(r/b)^a\). Consequently,
the conserved flux gives
\begin{equation}\label{eq:capacity-green-level}
 \Cap_p(\Omega_r)
 =\left(\frac ar\right)^{p-1}
       \int_{b=r}w^{p-1}\dd\sigma
 =a^{p-1}\om r^{n-p}.
\end{equation}
The normalized capacity used in \cite{BFM} is
\(\Cap_p(\Omega)/(a^{p-1}\om)\).
\begin{corollary}\label{cor:cap-curvature}
For every regular level, with \(H_+=\max\{H,0\}\),
\[
 \Cap_p(\Omega_r)
 \leq
 \left(\frac{a}{n-1}\right)^{p-1}
 |\Sigma_r|^{1/p}
 \left(\int_{b=r}H_+^p\dd\sigma\right)^{(p-1)/p}.
\]
\end{corollary}

\begin{proof}
Since \(\bs<1\), we may take \(\beta=1\) in
\eqref{eq:weighted-Minkowski}. H\"older's inequality gives
\[
 \frac{n-1}{r}\int_{b=r}w^p\dd\sigma
 \leq\int_{b=r}H_+w^{p-1}\dd\sigma
 \leq
 \left(\int_{b=r}H_+^p\dd\sigma\right)^{1/p}
 \left(\int_{b=r}w^p\dd\sigma\right)^{(p-1)/p}.
\]
Therefore
\[
 \int_{b=r}H_+^p\dd\sigma
 \geq\left(\frac{n-1}{r}\right)^p
       \int_{b=r}w^p\dd\sigma.
\]
A second application of H\"older's inequality and
the conserved flux give
\[
 \om r^{n-1}
 \leq|\Sigma_r|^{1/p}
       \left(\int_{b=r}w^p\dd\sigma\right)^{(p-1)/p}.
\]
Combining these inequalities with
\(\Cap_p(\Omega_r)=(a/r)^{p-1}\om r^{n-1}\)
proves the claim.
\end{proof}

The exterior derivative formula also measures the
geometric defect between scales.

\begin{theorem}[Annular defects and averaged umbilicity]
\label{thm:applications-main}
For \(\beta>\bs\),
\begin{equation}\label{eq:annular-main}
 R^{2-n}\int_{R<b<2R}\E\dd\mu
 \leq C(n,p,\beta)
       [A_\beta(R/2)-A_\beta(R)],
\end{equation}
where
\[
 C(n,p,\beta)
 =\frac{a(p-1)2^{n+a-2}}
        {\beta(1-2^{-a})}.
\]
The left-hand side tends to zero as \(R\to\infty\).
Moreover, if
\[
 1<p<2+\frac1{n-2},
\]
then
\[
 R\int_R^{2R}\frac1{|\Sigma_s|}
       \int_{b=s}|\mathring h|^2
                    \dd\sigma\dd s\longrightarrow0.
\]
In dimension three, this range is \(1<p<3\).
\end{theorem}

\begin{proof}
Integrating \eqref{eq:exterior} over \(t\in[R/2,R]\),
we obtain
\begin{align*}
 A_\beta(R/2)-A_\beta(R)
 &=\frac{\beta}{p-1}\int_{R/2}^R
       t^{a-1}\int_{b>t}
          b^{2-n-a}\E\dd\mu\dd t\\
 &\geq\frac{\beta}{p-1}
       \frac{R^a(1-2^{-a})}{a}(2R)^{2-n-a}
       \int_{R<b<2R}\E\dd\mu.
\end{align*}
This proves \eqref{eq:annular-main}.
Bounded monotonicity of \(A_\beta\) makes its dyadic
difference tend to zero.

Retaining the umbilicity term in \eqref{eq:E-def}
and applying coarea gives
\[
 R^{2-n}\int_R^{2R}\int_{b=s}
       w^{p+\beta-3}|\mathring h|^2
                    \dd\sigma\dd s\longrightarrow0.
\]
The choice \(\beta=3-p\) removes the gradient weight.
It lies in the strict parameter range precisely when
\[
 3-p>\frac{n-p}{n-1},
 \qquad\text{equivalently}\qquad
 p<2+\frac1{n-2}.
\]
Finally, the conserved flux and \(w\leq1\) give
\[
 \om s^{n-1}
 =\int_{b=s}w^{p-1}\dd\sigma
 \leq |\Sigma_s|
\]
for almost every \(s\). Hence
\[
 R\int_R^{2R}\frac1{|\Sigma_s|}
       \int_{b=s}|\mathring h|^2\dd\sigma\dd s
 \leq\om^{-1}R^{2-n}\int_R^{2R}
       \int_{b=s}|\mathring h|^2\dd\sigma\dd s,
\]
which tends to zero.
\end{proof}

These estimates control integral defects and averaged
umbilicity between scales. A metric approximation by cones,
a convergence rate, or uniqueness of a tangent cone requires
additional geometric arguments. The exponent \(3-p\) removes
the gradient weight and, in dimension three, also selects the
functional in scalar curvature monotonicity.

\subsection{Scalar curvature monotonicity in dimension three}
\label{subsec:cclt}

In dimension three, the choice \(\beta=3-p\) identifies the
Green-function functional studied by Munteanu--Wang \cite{MW}
and Chan--Chu--Lee--Tsang \cite{CCLT}. Write
\[
 a=\frac{3-p}{p-1},\qquad
 A(r)=A_{3-p}(r)=r^{-2}\int_{b=r}w^2\dd\sigma,
 \qquad M(r)=r^{-a}(A(r)-4\pi).
\]
If
\begin{equation}\label{eq:cclt-functional}
 \mathcal F(t)=\int_{G=t}|\nabla G|^2\dd\sigma,
 \qquad \Phi(t)=t^{-1}\mathcal F(t)-4\pi a^2t,
\end{equation}
then \(|\nabla G|=ar^{-a-1}w\) on \(b=r\) gives
\begin{equation}\label{eq:cclt-dictionary}
 \mathcal F(r^{-a})=a^2r^{-2a}A(r),\qquad
 \Phi(r^{-a})=a^2M(r).
\end{equation}
Thus increase of \(M\) in \(r\) corresponds to decrease of
\(\Phi\) in \(t\).

On the noncritical set, put
\begin{equation}\label{eq:cclt-Q}
 Q=|\mathring h|^2+\Ric(\nu,\nu)
   +|\nabla^\Sigma\log w|^2
   +\frac a2\left(H-\frac{2w}{b}\right)^2.
\end{equation}
Since \(3-p>\bs=(3-p)/2\), the exact formulas of
Section~\ref{sec:beta} apply for every \(1<p<3\).
Equation~\eqref{eq:E-def} gives
\(\mathcal E_{p,3-p}=wQ\); products involving \(Q\) below
are understood through this critical-set extension.

\begin{proposition}\label{prop:cclt-defect}
Under Assumption~\ref{ass:standing} with \(n=3\),
\begin{align}
 M'(r)&=a r^{-a-1}
       \int_{b<r}b^{-1}wQ\dd\mu,
 \label{eq:cclt-Mprime}\\
 A'(r)&=-a r^{a-1}
       \int_{b>r}b^{-a-1}wQ\dd\mu,
 \label{eq:cclt-Aprime}\\
 \Phi'(t)&=-a^2
       \int_{b<t^{-1/a}}b^{-1}wQ\dd\mu.
 \label{eq:cclt-Phiprime}
\end{align}
\end{proposition}

\begin{proof}
The first two identities are \eqref{eq:second-scaled}
and \eqref{eq:exterior}, with \(\beta/(p-1)=a\).
Differentiating \eqref{eq:cclt-dictionary} and using
\(dt/dr=-ar^{-a-1}\) gives
\[
 \Phi'(t)=-ar^{a+1}M'(r),
\]
which proves the third identity.
\end{proof}

In particular,
\begin{equation}\label{eq:cclt-Fprime}
 \mathcal F'(t)
 =t^{-1}\mathcal F(t)+4\pi a^2t
  -a^2t\int_{b<t^{-1/a}}b^{-1}wQ\dd\mu.
\end{equation}
Under \(\Ric\geq0\), the last integral is nonnegative.
This gives the CCLT comparison inequality in our curvature
setting and retains its separate geometric errors, for example
\begin{equation}\label{eq:cclt-small-defect}
 \int_{b<r}b^{-1}w|\mathring h|^2\dd\mu
 \leq\frac{r^{a+1}}a M'(r).
\end{equation}
CCLT work under nonnegative scalar curvature and
\(H_2(M;\mathbb Z)=0\). Their argument uses the Bochner formula,
an improved Kato inequality, Gauss--Bonnet, and regularization;
it includes first-variation identities and rigidity.
Their original theorem treats \(1<p\leq2\), with an additional
smoothness assumption for \(2<p<3\). The later analysis in
\cite[Appendix~A]{BFMmass} also addresses critical points for
the related scalar curvature monotonicities in this latter range.
Our formulas provide a nonnegative volume remainder under
\(\Ric\geq0\), within the full Colding-type family. They do not
establish the scalar curvature theorem in its broader curvature
class. Under our assumptions, \(M'(r_0)=0\) gives global Euclidean
rigidity, whereas \(A'(r_0)=0\) at a regular level gives exterior
cone rigidity, by Theorem~\ref{thm:beta-rigidity}.

The scalar curvature structure itself can be read from the
same local identity. The following computation requires no
curvature sign.

\begin{proposition}\label{prop:scalar-bridge}
Let \(n=3\), and let \(b>0\) be smooth with \(w>0\) on a band
foliated by compact boundaryless levels. Suppose
\(\Delta_pb=2w^p/b\), and define \(A,M\) as above. Set
\begin{equation}\label{eq:cclt-P}
 \mathcal P(r)=\int_{b=r}\left[
   \frac{\operatorname{Scal}}2+\frac{|\mathring h|^2}{2}
   +|\nabla^\Sigma\log w|^2
   +\frac{2a+1}{4}\left(H-\frac{2w}{r}\right)^2
   \right]\dd\sigma.
\end{equation}
Then
\begin{equation}\label{eq:scalar-bridge}
 \left(r^{a+2}M'(r)\right)'
 =a\left[\mathcal P(r)+4\pi-2\pi\chi(\Sigma_r)\right].
\end{equation}
\end{proposition}

\begin{proof}
The local divergence identity \eqref{eq:master-div} is valid
without a curvature sign. Integrating it over a band, using
its flux \((p-1)r^{1-a}A'(r)\), and applying coarea gives
\[
 \left(r^{1-a}A'(r)\right)'
 =a r^{-a-1}\int_{b=r}Q\dd\sigma.
\]
Multiplication by \(r^{a+1}\) yields
\begin{equation}\label{eq:cclt-ODE1}
 r^2A''+(1-a)rA'=a\int_{b=r}Q\dd\sigma.
\end{equation}
Also, \eqref{eq:A-H} is a local first-variation identity, so
\begin{equation}\label{eq:cclt-ODE2}
 A'=-ar^{-2}\int_{b=r}w\left(H-\frac{2w}{r}\right)\dd\sigma.
\end{equation}
The Gauss equation and an expansion of the square give
\begin{align*}
 \Ric(\nu,\nu)
 &=\frac{\operatorname{Scal}}2-K_\Sigma
     +\frac{H^2}{4}-\frac{|\mathring h|^2}{2},\\
 \frac{H^2}{4}+\frac a2\left(H-\frac{2w}{r}\right)^2
 &=\frac{2a+1}{4}\left(H-\frac{2w}{r}\right)^2
     +\frac wr\left(H-\frac{2w}{r}\right)+\frac{w^2}{r^2}.
\end{align*}
Integration, Gauss--Bonnet, and \eqref{eq:cclt-ODE2} imply
\[
 \int_{b=r}Q\dd\sigma
 =\mathcal P(r)-2\pi\chi(\Sigma_r)-\frac r a A'(r)+A(r).
\]
Substitution into \eqref{eq:cclt-ODE1} gives
\[
 r^2A''+(2-a)rA'-a(A-4\pi)
 =a[\mathcal P(r)+4\pi-2\pi\chi(\Sigma_r)].
\]
This is \eqref{eq:scalar-bridge}, because
\(r^{a+2}M'=r^2A'-ar(A-4\pi)\).
\end{proof}

For \(\operatorname{Scal}\geq0\) and connected levels,
\(\mathcal P\geq0\) and \(\chi(\Sigma_r)\leq2\).
This explains the role of level-set topology in scalar
curvature monotonicity. The passage across critical levels
requires the additional analysis in the cited works.

\subsection{Connections with nonlinear potential theory}
\label{subsec:nonlinear-potential}

Mari--Rigoli--Setti \cite{MRS} provide a direct link
in the original metric. In the Euclidean model, their fake
distance is \(b\) and their radial weight exponent is
\(2-n-a\). Their equation~(2.28) becomes, on \(w>0\),
\[
 \frac12 b^{n+a-2}
 \Div\bigl(b^{2-n-a}w^{p-2}\A\nabla w^2\bigr)
 =\mathcal E_{p,2}.
\]
Thus our \(\beta=2\) divergence identity is already
contained in \cite[Proposition~2.15]{MRS}. Their auxiliary
power is proportional to \(b^{1+1/\sqrt{p-1}}\), and
their anisotropic Hessian norm differs, for \(p\ne2\),
from the norm \(|\cdot|_{\A}\) used here with
\(v=b^{p/(p-1)}\). These choices give different square
decompositions of the same scalar identity at \(\beta=2\).
Their normalized surface average in Lemma~5.3 equals
\(I_u/\om\); for \(u=F\), their surface and volume
averages are \(A_p/\om\) and \(nV_p/\om\).

Fogagnolo--Mazzieri--Pinamonti \cite{FMP} proved monotonicity
and rigidity for exterior \(p\)-capacitary potentials of smooth
convex Euclidean domains. Their conformal calculation uses a
generalized \(p\)-Bochner formula and a decomposition of the
Hessian; see \cite[Lemma~4.2 and Proposition~4.4]{FMP}.
Benatti--Fogagnolo--Mazzieri extended this theory to
\(p\)-nonparabolic manifolds with \(\Ric\geq0\), including
critical points, and proved exact derivative and convexity
identities with exterior cone rigidity
\cite[Theorems~3.1--3.2]{BFM}.
For \(1<p\leq2\), \cite[Theorem~3.1]{BPP} uses the same
three geometric squares as \eqref{eq:E-def}, with the Ricci
term treated separately and their parameter equal to our
\(\beta\). They also establish curvature-varifold regularity
of almost every level.

The relation with our family is explicit. Choose a small
regular value \(r_0\), and let \(u=(r_0/b)^a\) be the exterior
capacitary potential of \(\Omega_{r_0}\). In the notation of
\cite{BFM},
\[
 \mathcal F_p^\eta(t)
 =t^{\eta(n-1)(p-1)/(n-p)}
       \int_{u=1/t}|\nabla u|^{(\eta+1)(p-1)}\dd\sigma.
\]
For \(r\geq r_0\), setting \(t=(r/r_0)^a\) and using
\(|\nabla u|=ar_0^ab^{-a-1}w\) gives
\begin{equation}\label{eq:BFM-dictionary}
 \mathcal F_p^\eta((r/r_0)^a)
 =a^{(\eta+1)(p-1)}r_0^{n-p-\eta(p-1)}
       A_{\eta(p-1)}(r).
\end{equation}
Thus \(\beta=\eta(p-1)\), and the strict parameter ranges
agree. For these Green-level exteriors, the third formula
is their convexity identity after reparametrization, and the
exterior formula is their first-derivative identity. The other
two formulas follow by weighted integration with the normalized
pole limits, as computed in Sections~\ref{sec:fixed-formulas}
and \ref{sec:beta}.

The local defects also correspond. For
\(\widetilde g=(r_0/b)^2g\) and
\(\varphi=(n-2)\log(b/r_0)\), the conformal Hessian formula
and \eqref{eq:B-T} give
\begin{equation}\label{eq:potential-B}
 B=\frac{q^{p-1}}{n-2}b^2w^{p-2}
                 \widetilde\nabla^{\,2}\varphi,
 \qquad \Delta_{p,\widetilde g}\varphi=0.
\end{equation}
\rev{Transforming both the flux and the volume form in
\cite[Lemma~3.6]{BFM} back to \(g\) gives
\eqref{eq:master-div}, up to a fixed positive normalization
factor, and hence also \eqref{eq:master-product}.}
These correspondences show that the weighted
original-metric and conformal calculations share the same
geometric defect structure. Our formulation organizes this
structure around \(v=b^q\) and \(\LG(Gw^\beta)\): the first
linearization produces \(B\), and the second removes the
zero-order term and matches the conserved flux. This gives
a single derivation of Colding's three surface and volume
formulas, including the normalized pole terms, equality
cases, and annular defect estimates. The sharp gradient
bound, already proved in \cite{MRS}, can also be recovered
from the supremum monotonicity in \cite{BFM} and the
uniform pole limit.

\subsection{Relations with mass and the Penrose inequality}

Still in dimension three, let \(u=(r_0/b)^a\) and
\(\psi=-(p-1)\log u\). Then
\[
 \Delta_p\psi=|\nabla\psi|^p,
 \qquad \frac{|\nabla\psi|}{3-p}=\frac wb.
\]
The formal limit as \(p\downarrow1\) is the level-set
equation for inverse mean curvature flow
\cite{Moser,KN,MRS}.
Rigorous convergence of the potentials, their gradients,
and the associated monotonicity formulas is proved under the
assumptions in \cite{BPP}; no \(p\downarrow1\) limit is taken
in the present paper.
With the normalized capacity
\(c_p(\Omega)=a^{1-p}\Cap_p(\Omega)/(4\pi)\),
\eqref{eq:capacity-green-level} gives
\(c_p(\Omega_r)=r^{3-p}\).
Consequently, the \(p\)-Hawking mass of
\cite[Section~2.2]{BFMmass} becomes, on a regular Green level,
\begin{equation}\label{eq:green-hawking}
 m_H^{(p)}(\Sigma_r)
 =\frac r{8\pi}\left[
   4\pi+A(r)-\frac1r\int_{b=r}Hw\dd\sigma\right].
\end{equation}
This uses the capacity factor \(c_p(\Omega_r)^{1/(3-p)}=r\)
and the relation \(|\nabla\psi|/(3-p)=w/r\) on the level.
By \eqref{eq:cclt-ODE2},
\[
 \frac1r\int_{b=r}Hw\dd\sigma=2A(r)-\frac r a A'(r),
\]
so the derivative--mass correspondence reads
\begin{equation}\label{eq:mass-derivative-dictionary}
 m_H^{(p)}(\Sigma_r)
 =\frac r{8\pi}\left[4\pi-A(r)+\frac r a A'(r)\right]
 =\frac{r^{a+2}}{8\pi a}M'(r).
\end{equation}
This is the correspondence in \cite[Appendix~A]{BFMmass}
written in the Green radius. Under our Ricci assumptions,
\eqref{eq:cclt-Mprime} further yields
\begin{equation}\label{eq:mass-volume-defect}
 m_H^{(p)}(\Sigma_r)
 =\frac r{8\pi}\int_{b<r}b^{-1}wQ\dd\mu.
\end{equation}
On a smooth band, the local scalar identity
\eqref{eq:scalar-bridge} gives
\begin{equation}\label{eq:mass-scalar-derivative}
 \frac d{dr}m_H^{(p)}(\Sigma_r)
 =\frac1{8\pi}[\mathcal P(r)+4\pi-2\pi\chi(\Sigma_r)].
\end{equation}
Its right-hand side is nonnegative for nonnegative scalar
curvature and connected levels.

These are correspondences with potential-theoretic mass
monotonicity. The global mass and Penrose inequalities in
\cite{BFMmass,BFMPenrose,HMT,XXZ} use boundary, topological,
and asymptotic assumptions beyond the present framework.
Under Assumption~\ref{ass:standing}, asymptotic flatness
would give asymptotic volume ratio one, so Bishop--Gromov
rigidity would force \(M\) to be Euclidean. Nontrivial
mass problems therefore require a broader curvature setting;
positive-mass spatial Schwarzschild metrics, for example,
have zero scalar curvature and negative radial Ricci curvature.

\medskip
	\noindent\textbf{Disclosure on AI assistance.}
The main ideas and research motivation of this work
originate with the author. ChatGPT assisted with
parts of the proofs and with the revision of the
manuscript. The author reviewed and verified the
mathematical arguments and the final text, and
takes full responsibility for the content and
correctness of the paper.

\begin{flushleft}
Yu-Zhao Wang\\
School of Mathematics and Statistics, Shanxi University\\
Taiyuan 030006, Shanxi, China\\
\emph{Email}: \href{mailto:wangyuzhao@sxu.edu.cn}{wangyuzhao@sxu.edu.cn}
\end{flushleft}

\end{document}